\RequirePackage{fix-cm}
\documentclass[smallextended]{svjour3}       
\smartqed  
\usepackage{graphicx}
\usepackage{amsmath,amsfonts}
\usepackage{amssymb}
\usepackage{latexsym}
\usepackage{mathrsfs}
\usepackage{colortbl}
\usepackage{amscd}
\usepackage{tikz-cd}
\usepackage{comment}
\usepackage{url}
\usepackage[subrefformat=parens]{subcaption}
\usepackage{mathtools}
\usepackage[section]{placeins}

\allowdisplaybreaks[3]

  \journalname{}

\def\diam{\mathop{\mathrm{diam}}}
\def\diag{\mathop{\mathrm{diag}}}

\def\div{\mathop{\mathrm{div}}}

\def\<{\mathop{\textless}}
\def\>{\mathop{\textgreater}}

\def\card{\mathop{\rm{card}}}

\spnewtheorem{thr}{Theorem}{\bf}{\it}
\spnewtheorem{defi}{Definition}{\bf}{\it}
\spnewtheorem{lem}{Lemma}{\bf}{\it}
\spnewtheorem{coro}{Corollary}{\bf}{\it}
\spnewtheorem{assume}{Assumption}{\bf}{\it}
\spnewtheorem{ex}{Example}{\bf}{\it}
\spnewtheorem{Case}{Case}{\bf}{\it}
\spnewtheorem*{pf*}{Proof}{\bf}{\rm}
\spnewtheorem*{rem*}{Remark:}{\it}{\it}
\spnewtheorem*{ex*}{Example:}{\it}{\it}
\spnewtheorem{Cond}{Condition}{\bf}{\it}
\spnewtheorem{rem}{Remark}{\it}{\it}

\spnewtheorem*{lem1*}{Lemma 1}{\bf}{\rm}
\spnewtheorem*{appendix1*}{Appendix A}{\bf}{\rm}

\newcounter{sone}
\newcounter{stwo}
\newcounter{sthree}
\newcounter{sfour}
\newcounter{sfive}
\newcounter{ssix}
\newcounter{lone}
\newcounter{ltwo}
\newcounter{lthree}
\newcounter{lfour}
\newcounter{lfive}
\newcounter{lsix}
\makeatletter
    
    \@addtoreset{equation}{section}

\begin{document}

\title{Hybrid weakly over-penalised symmetric interior penalty method on anisotropic meshes
}

\titlerunning{Hybrid WOPSIP methods on anisotropic meshes}        

\author{Hiroki Ishizaka  
}


\institute{Hiroki Ishizaka \at
              Team FEM, Matsuyama, Japan \\
              \email{h.ishizaka005@gmail.com}           
}

\date{Received: date / Accepted: date}

\maketitle

\begin{abstract}
In this study, we investigate a hybrid-type anisotropic weakly over-penalised symmetric interior penalty method for the Poisson equation on convex domains. Compared with the well-known hybrid discontinuous Galerkin methods, our approach is simple and easy to implement. Our primary contributions are the proposal of a new scheme and the demonstration of a proof for the consistency term, which allows us to estimate the anisotropic consistency error. The key idea of the proof is to apply the relation between the Raviart--Thomas finite element space and a discontinuous space. In numerical experiments, we compare the calculation results for standard and anisotropic mesh partitions.

\end{abstract}

\section{Introduction} \label{intro}
In this study, we investigate a hybrid weakly over-penalised symmetric interior penalty (HWOPSIP) method for Poisson equations on anisotropic meshes. Brenner et al. were the first to propose the original weakly over-penalised symmetric interior penalty (WOPSIP) method \cite{BreOweSun08}. Several studies have considered similar techniques since \cite{BarBre14,Bre15,BreOweSun12}. WOPSIP methods have several advantages over the standard symmetric interior penalty (SIP) discontinuous Galerkin \cite{PieErn12,Riv08}, and the classic Crouzeix--Raviart (CR) nonconforming finite element \cite{CroRav73} methods. One advantage is that the WOPSIP method is stable for any penalty parameter. The WOPSIP method is similar to the classical CR finite element method and thus has similar features, such as inf-sup stability on anisotropic meshes (see \cite[Lemma 7]{Ish24}). Another advantage of the WOPSIP method is that error analysis is performed on more general meshes (\cite{BreOweSun12,BreOweSun08,Bre15,BarBre14}) than conformal meshes. Applying nonconforming meshes can be difficult in the CR nonconforming finite element method. Furthermore, a consistency term is not included in the WOPSIP scheme. This may be advantageous when using anisotropic meshes; that is, the use of WOPSIP methods may prevent worsening interpolation errors because of the trace inequality (Lemma \ref{lem=trace}).

Anisotropic meshes are effective for problems in which the solution exhibits anisotropic behaviour in some directions of the domain. Anisotropic meshes can be divided into two types: those that include elements with large aspect ratios and those that do not satisfy the shape-regularity condition but satisfy the semi-regular condition ({Assumption} \ref{neogeo=assume}); and those that include elements with large aspect ratios and whose partitions satisfy the shape-regularity condition. The former includes  graded meshes, whereas the latter includes Shishkin meshes. In several studies, isotropic mesh partitions that satisfy the shape-regularity condition are used for analysis, that is, triangles or tetrahedra cannot be overly flat in a shape-regular family of triangulations. In studies on anisotropic finite element methods, it is crucial to deduce anisotropic interpolation error estimates and sharper error estimates than those in standard analysis. In \cite{Ish24}, we studied anisotropic error analysis for the Stokes using the WOPSIP method.

The hybrid discontinuous Galerkin (HdG) method proposed by Cockburn et al. \cite{Cocetal09} has received considerable attention in recent years. In SIP methods, the numerical traces of functions are used as functions of the interfaces. On the other hand, the main idea of the HdG methods is to introduce additional numerical traces, which are single-valued functions on an interface. In this study, we apply this idea to the WOPSIP method and perform error estimations for both the energy and $L^2$ norms under a semi-regular condition, which is equivalent to the maximum angle condition. The proposed scheme \eqref{hwop=4} is known as the HWOPSIP method. The error between the exact and HWOPSIP finite element approximation solutions with the energy norm is divided into two parts. One part is the optimal approximation error in discontinuous finite element spaces and the other is a consistency error term. For the former, first-order CR interpolation errors (Theorem \ref{DGCE=thr3}) are used. However, estimating the consistency error term on anisotropic meshes is challenging. On anisotropic meshes, analyses using trace inequalities do not yield optimal estimates. This difficulty can be overcome using the relation between the lowest-order Raviart--Thomas (RT) finite element interpolation and discontinuous space; see Lemma \ref{rel=CRRT}. This relation derives an optimal error estimate for the consistency error (Lemma \ref{asy=con}).

The remainder of this paper is organised as follows. In Section 2, we introduce the Poisson equation with the Dirichlet boundary condition and its weak formulation and present the proposed HWOPSIP approximation for the continuous problem. In Section 3, we introduce the finite element settings of the HWOPSIP method. In Section 4, the discrete Poincar\'e inequality is presented. In Section 5, we discuss the stability and error estimates of the HWOPSIP approximate problem. In Section 6, we provide the results of numerical evaluations are ptovided. Finally, in Section 7, we conclude the paper by noting some limitations of the present study, suggesting possible directions for future research. 

Throughout this paper, we denote by $c$ a constant independent of $h$ (defined later) and the angles and aspect ratios of simplices, unless specified otherwise {all constants $c$ are bounded if the maximum angle is bounded}. These values vary in different contexts. The notation $\mathbb{R}_+$ denotes the set of positive real numbers.

\section{HWOPSIP method for the Poisson equation} \label{sec=HyWOPSIP}
This section presents an analysis of the HWOPSIP method for the Poisson equations on anisotropic meshes.

\subsection{Continuous problems}
Let $\Omega \subset \mathbb{R}^d$, $d \in \{ 2 , 3 \}$ be a bounded polyhedral domain. Furthermore, we assume that $\Omega$ is convex if necessary. The Poisson problem determines $u: \Omega \to \mathbb{R}$ such that
\begin{align}
\displaystyle
- \varDelta u  = f \quad \text{in $\Omega$}, \quad u = 0 \quad \text{on $\partial \Omega$}, \label{poisson_eq}
\end{align}
where $f \in L^2(\Omega)$ is a given function. The variational formulation of the Poisson equations \eqref{poisson_eq} is as follows: $u \in H_0^1(\Omega)$ is determined such that
\begin{align}
\displaystyle
a(u,\varphi) := \int_{\Omega} \nabla u \cdot \nabla \varphi dx &= \int_{\Omega} f  \varphi dx \quad \forall \varphi \in H_0^1(\Omega). \label{poisson_weak}
\end{align}
By the Lax--Milgram lemma, there exists a unique solution $u \in H_0^1(\Omega)$ for any $f \in L^2(\Omega)$ and it holds that
\begin{align*}
\displaystyle
| u |_{H^1(\Omega)} \leq C_P(\Omega) \| f \|,
\end{align*}
where $C_P(\Omega)$ is the Poincar$\rm{\acute{e}}$ constant depending on $\Omega$. 
Furthermore, if $\Omega$ is convex, then $u \in H^2(\Omega)$ and 
\begin{align}
\displaystyle
| u |_{H^2(\Omega)} \leq \| \varDelta u \|. \label{poisson=reg}
\end{align}
The proof can be found in \cite[Theorem 3.1.1.2, Theorem 3.2.1.2]{Gri11}.

\subsection{Meshes and mesh faces} \label{regularmesh0}
Let $\mathbb{T}_h = \{ T \}$ be a simplicial mesh of $\overline{\Omega}$ composed of closed $d$-simplices, such as $\overline{\Omega} = \bigcup_{T \in \mathbb{T}_h} T$, with $h := \max_{T \in \mathbb{T}_h} h_{T}$, where $ h_{T} := \diam(T)$. For simplicity, we assume that $\mathbb{T}_h$ is conformal, that is, $\mathbb{T}_h$ is a simplicial mesh of $\overline{\Omega}$ without hanging nodes.

Let $\mathcal{F}_h^i$ be the set of interior faces, and $\mathcal{F}_h^{\partial}$ be the set of faces on boundary $\partial \Omega$. We set $\mathcal{F}_h := \mathcal{F}_h^i \cup \mathcal{F}_h^{\partial}$. For any $F \in \mathcal{F}_h$, we define the unit normal $n_F$ to $F$ as follows: (\roman{sone}) If  $F \in \mathcal{F}_h^i$ with $F = T_1 \cap T_2$, $T_1,T_2 \in \mathbb{T}_h$, let $n_1$ and $n_2$ be the outward unit normals of $T_1$ and $T_2$, respectively. Then, $n_F$ is either $\{ n_1 , n_2\}$; (\roman{stwo}) If $F \in \mathcal{F}_h^{\partial}$, $n_F$ is the unit outwards normal $n$ to $\partial \Omega$. Furthermore, for a simplex $T \subset \mathbb{R}^d$, let $\mathcal{F}_{T}$ be the collection of the faces of $T$.

For any $k \in \mathbb{N}_0$, let $\mathbb{P}^k(T)$ and $\mathbb{P}^k(F)$ be spaces of polynomials with a degree of at most $k$ in $T$ and $F$, respectively. For any $F \in \mathcal{F}_h$, we define the $L^2$-projection $\Pi_F^{0}: L^2(F) \to \mathbb{P}^{0}(F)$ as
\begin{align*}
\displaystyle
\int_F (\Pi_F^{0} \varphi - \varphi)   ds = 0 \quad \forall \varphi \in L^2(F).
\end{align*}

\subsection{Penalty parameter}
We consider an appropriate parameter for the penalty term in HdG methods on anisotropic meshes. The following trace inequality on anisotropic meshes is significant in this study: 

\begin{lem}[Trace inequality] \label{lem=trace}
Let $T \subset \mathbb{R}^d$ be a simplex.  Then, there exists a positive constant $c$ such that for any $v \in H^{1}(T)^d$, $F \in \mathcal{F}_{T}$, and $h$
\begin{align}
\displaystyle
\| v \|_{L^2(F)^d}
\leq c \ell_{T,F}^{- \frac{1}{2}} \left( \| v \|_{L^2(T)^d} + h_{T}^{\frac{1}{2}}  \| v \|_{L^2(T)^d}^{\frac{1}{2}} | v |_{H^1(T)^d}^{\frac{1}{2}} \right), \label{trace=vec}
\end{align}
where $\ell_{T,F} := \frac{d! |T|_d}{|F|_{d-1}}$ denotes the distance from the vertex of $T$ opposite to $F$ to the face.
\end{lem}

\begin{pf*}
The proof can be found in \cite[Lemma 1]{Ish24}.
\qed
\end{pf*}

{
We define a broken (piecewise) Hilbert space as
\begin{align*}
\displaystyle
H^{1}(\mathbb{T}_h) := \{ v \in L^2(\Omega): \ v|_{T} \in H^{1}(T) \quad \forall T \in \mathbb{T}_h  \}
\end{align*}
with a norm
\begin{align*}
\displaystyle
| v |_{H^{1}(\mathbb{T}_h)} &:= \left( \sum_{T \in \mathbb{T}_h} | v |^2_{H^{1}(T) } \right)^{\frac{1}{2}}.
\end{align*}
}
Let $\beta$ be nonnegative real numbers chosen latter on. Let $v \in H^{1}(\mathbb{T}_h)^d$ and $\varphi \in L^{2}(F)$. Using the trace and the H\"older inequalities gives the following estimate for the term $\int_F | v \cdot n_F  \Pi_F \varphi | ds$,
\begin{align*}
\displaystyle
\int_F | v \cdot n_F  \Pi_F \varphi | ds
&\leq c \ell_{T,F}^{- \frac{1}{2}} \left( \| v \|_{L^2(T)^d} + h_{T}^{\frac{1}{2}}  \| v \|_{L^2(T)^d}^{\frac{1}{2}} | v |_{H^1(T)^d}^{\frac{1}{2}} \right) \| \Pi_F \varphi \|_{L^2(F)} \\
&\hspace{-1cm} = c h^{\beta} \left( \| v \|_{L^2(T)^d} + h_{T}^{\frac{1}{2}}  \| v \|_{L^2(T)^d}^{\frac{1}{2}} | v |_{H^1(T)^d}^{\frac{1}{2}} \right) h^{-\beta } \ell_{T,F}^{- \frac{1}{2}} \| \Pi_F \varphi \|_{L^2(F)}.
\end{align*}
The associated penalty parameter $\kappa_{F(\beta)}$ becomes 
\begin{align*}
\displaystyle
\kappa_{F(\beta)} := h^{- 2 \beta} \ell_{T,F}^{- 1}.
\end{align*}
{The penalty parameter is used in the discrete Poincar\'e inequality (Section \ref{disPoi}) when $\beta = 0$, and in the HWOPSIP method (Section \ref{new=sch}) when $\beta = 1$.}

\subsection{Finite element spaces}
For $s \in \mathbb{N}_0$, we define the standard discontinuous finite-element space as
\begin{align*}
\displaystyle
P_{dc,h}^{s} &:= \left\{ p_h \in L^2(\Omega); \ p_h|_{T} \in \mathbb{P}^{s}({T}) \quad \forall T \in \mathbb{T}_h  \right\}.
\end{align*}

Let $Ne$ be the number of elements included in mesh $\mathbb{T}_h$. Thus, we have $\mathbb{T}_h = \{ T_j\}_{j=1}^{Ne}$. Using the barycentric coordinates $ \{{\lambda}_{T_j,i} \}_{i=1}^{d+1}: \mathbb{R}^d \to \mathbb{R}$, the nodal basis functions for the CR finite element method are defined as
\begin{align}
\displaystyle
{\theta}^{CR}_{T_j,i}({x}) := d \left( \frac{1}{d} - {\lambda}_{T_j,i} ({x}) \right) \quad \forall i \in \{ 1, \ldots ,d+1 \}. \label{CR2}
\end{align}
For $j \in \{1, \ldots ,Ne \}$ and $i \in \{ 1, \ldots , d+1\}$, we define the function $\phi_{j(i)}$ as
\begin{align}
\displaystyle
\phi_{j(i)}(x) :=
\begin{cases}
\theta_{T_j,i}^{CR}(x), \quad \text{$x \in T_j$}, \\
0, \quad \text{$x \notin T_j$}.
\end{cases} \label{CR5}
\end{align}
We define a discontinuous finite-element space as
\begin{align}
\displaystyle
V_{dc,h}^{CR} &:= \left\{ \sum_{j=1}^{Ne} \sum_{i=1}^{d+1} c_{j(i)} \phi_{j(i)}; \  c_{j(i)} \in \mathbb{R}, \ \forall i,j \right\} \subset P_{dc,h}^1. \label{CR6}
\end{align}
We define the function space as
\begin{align*}
\displaystyle
L_0^2(\mathcal{F}_h) := \left \{ \lambda \in L^2(\mathcal{F}_h): \ \lambda = 0 \ \text{on $\partial \Omega$}  \right \}.
\end{align*}
Note that $\lambda \in L^2(\mathcal{F}_h)$ is a single-valued function on $F \in \mathcal{F}_h$. We set $V^H(h) := H^1(\mathbb{T}_h) \times L^2(\mathcal{F}_h)$ with the following norm:
\begin{align*}
\displaystyle
| v^H |_{hwop(\beta)} &:= \left( |v|^2_{H^1(\mathbb{T}_h)} + |v^H |^2_{jhwop(\beta)} \right)^{\frac{1}{2}},\\
 |v^H|_{jhwop(\beta)} &:= \left( \sum_{T \in \mathbb{T}_h} \sum_{F \in \mathcal{F}_{T}} \kappa_{F(\beta)} \| \Pi_F^0 (v - \lambda) \|_{L^2(F)}^2 \right)^{\frac{1}{2}}
\end{align*}
for any $v^H := (v, \lambda) \in V^{H}(h)$. We set $V^H_0(h) := H^1(\mathbb{T}_h) \times L_0^2(\mathcal{F}_h)$, $X^H := H^1(\Omega) \times L^2(\mathcal{F}_h)$, and $X^H_0 := H^1(\Omega) \times L^2_0(\mathcal{F}_h)$. Furthermore, we define spaces as
\begin{align*}
\displaystyle
V_{dc,h}(\mathcal{F}_h) &:= \left\{ \lambda_h \in L^2(\mathcal{F}_h): \ \lambda_h|_{F} \in \mathbb{P}^0(F) \ \forall F \in \mathcal{F}_h  \right\} \subset L^2(\mathcal{F}_h), \\ 
V_{dc,h,0}(\mathcal{F}_h) &:= V_{dc,h}(\mathcal{F}_h) \cap L_0^2(\mathcal{F}_h) \subset L_0^2(\mathcal{F}_h), \\
V_{dc,h,0}^{H} &:= V_{dc,h}^{CR} \times V_{dc,h,0}(\mathcal{F}_h) \subset  V^H_0(h), \\
V_{dc,h}^{H} &:= V_{dc,h}^{CR} \times V_{dc,h}(\mathcal{F}_h) \subset  V^H(h).
\end{align*}

\subsection{New scheme} \label{new=sch}
We consider the HWOPSIP method for Poisson equation \eqref{poisson_eq} as follows: We aim to determine $u_h^H := (u_h^{hwop}, \lambda_h^{hwop}) \in V_{dc,h,0}^{H}$ such that
\begin{align}
\displaystyle
a_h^{hwop}(u_h^{H} , \varphi_h^H) = \ell_h(\varphi_h^H) \quad \forall \varphi_h^H := (\varphi_h,\mu_h) \in V_{dc,h,0}^{H}, \label{hwop=4}
\end{align}
where $a_h^{hwop}: (V_{dc,h}^{H} + X^H ) \times (V_{dc,h}^{H} + X^H ) \to \mathbb{R}$ and $\ell_h: V_{dc,h}^{H} + X^H  \to \mathbb{R}$ are defined as
\begin{align*}
\displaystyle
a_h^{hwop}(v^H,\varphi^H) &:= \int_{\Omega} \nabla_h v \cdot \nabla_h \varphi dx \\
&\quad + \sum_{T \in \mathbb{T}_h} \sum_{F \in {\mathcal{F}_T}} \kappa_{F(1)}  \int_F \Pi_F^{0} ( v - \lambda )\Pi_F^{0} ( \varphi - \mu ) ds, \\
\ell_h(\varphi^H) &:= \int_{\Omega} f \varphi dx
\end{align*}
for all $v^H :=(v,\lambda) \in V_{dc,h}^{H} + X^H $ and $\varphi^H :=(\varphi,\mu) \in V_{dc,h}^{H} + X^H $. Here, we define a broken gradient operator as follows. For $\varphi \in H^1(\mathbb{T}_h)$, the broken gradient $\nabla_h:H^1(\mathbb{T}_h) \to L^2(\Omega)^{d}$ is defined as
\begin{align*}
\displaystyle
(\nabla_h \varphi)|_T &:= \nabla (\varphi|_T) \quad \forall T \in \mathbb{T}_h.
\end{align*}
Using H\"older's inequality, we obtain
\begin{align}
\displaystyle
|a_h^{hwop}(v^H,\varphi^H)|
&\leq c |v^H|_{hwop(1)} |\varphi^H|_{hwop(1)} \label{hwop=5}
\end{align}
for all $v^H  \in V_{dc,h}^{H} + X^H $ and $\varphi^H  \in V_{dc,h}^{H} + X^H $. For any $\varphi_h^H := (\varphi_h,\mu_h) \in V_{dc,h,0}^{H}$, it holds that
\begin{align}
\displaystyle
a_h^{hwop}(\varphi_h^{H} , \varphi_h^H) = | \varphi_h^H |_{hwop(1)}^2. \label{hwop=6}
\end{align}

\section{Anisotropic interpolation error estimates} \label{Ani=int}

\subsection{Reference elements} \label{reference}
We now define the reference elements $\widehat{T} \subset \mathbb{R}^d$.

\subsubsection*{Two-dimensional case} \label{reference2d}
Let $\widehat{T} \subset \mathbb{R}^2$ be a reference triangle with vertices $\hat{p}_1 := (0,0)^T$, $\hat{p}_2 := (1,0)^T$, and $\hat{p}_3 := (0,1)^T$. 

\subsubsection*{Three-dimensional case} \label{reference3d}
In the three-dimensional case, we consider the following two cases: (\roman{sone}) and (\roman{stwo}); see Condition \ref{cond2}.

Let $\widehat{T}_1$ and $\widehat{T}_2$ be reference tetrahedra with the following vertices:
\begin{description}
   \item[(\roman{sone})] $\widehat{T}_1$ has the vertices $\hat{p}_1 := (0,0,0)^T$, $\hat{p}_2 := (1,0,0)^T$, $\hat{p}_3 := (0,1,0)^T$, and $\hat{p}_4 := (0,0,1)^T$;
 \item[(\roman{stwo})] $\widehat{T}_2$ has the vertices $\hat{p}_1 := (0,0,0)^T$, $\hat{p}_2 := (1,0,0)^T$, $\hat{p}_3 := (1,1,0)^T$, and $\hat{p}_4 := (0,0,1)^T$.
\end{description}
Therefore, we set $\widehat{T} \in \{ \widehat{T}_1 , \widehat{T}_2 \}$. 
Case (\roman{sone}) is called the \textit{ regular vertex property}; see \cite{AcoDur99}.

\subsection{Standard affine mapping}
In standard interpolation theory, an affine mapping that connects the reference element to the mesh element is introduced. However, interpolation errors may be overestimated in anisotropic meshes. 

Let $\widehat{T} \subset \mathbb{R}^d$ and  $T \subset \mathbb{R}^d$ be a reference element and a simplex, respectively. Let these two elements be affine equivalent. The transformation $\Psi_{T}$ takes the form
\begin{align*}
\displaystyle
&\Psi_{T}: \widehat{T} \ni \hat{x}  \mapsto \Psi_{T}(\hat{x}) := {B}_{T} \hat{x} + b_{T} \in T, 
\end{align*}
where ${B}_{T} \in \mathbb{R}^{d \times d}$ is an invertible matrix and $b_{T} \in \mathbb{R}^{d}$. In classical theory (e.g. \cite[Section 1.5.1]{ErnGue04}), the quantity $ \| {B}_{T} \|_2 \| {B}_{T}^{-1} \|_2$, which is called the \textit{Euclidean condition number} of ${B}_{T}$, may be included in standard interpolation error estimates. 
\begin{rem}
Using a standard argument ( see \cite[Lemma 1.100]{ErnGue04}), we have
\begin{align*}
\displaystyle
\|  {B}_{T} \|_2 \| {B}_{T}^{-1} \|_2 \leq c \frac{h_{T}}{\rho_{T}},
\end{align*}
where $\|  {B}_{T} \|_2$ is the operator norm of ${B}_T$, and $\rho_T$ is the diameter of the largest ball that can be inscribed in $T$. For a shape-regular family of meshes, that is, there exists a constant $\gamma_1 \> 0$ such that
\begin{align}
\displaystyle
\rho_{T} \geq \gamma_1 h_{T} \quad \forall \mathbb{T}_h \in \{ \mathbb{T}_h \}, \quad \forall T \in \mathbb{T}_h, \label{geo1}
\end{align}
quantity $ \| {B}_{T} \|_2 \| {B}_{T}^{-1} \|_2$ is bounded. However, when the shape-regularity condition is violated (i.e. the simplex becomes as flat as $h_{T} \to 0$), and the quantity may diverge.
\end{rem}

\color{black}
\subsection{Two-step affine mapping} \label{two=step}
We introduce a new strategy proposed in \cite[Section 2]{IshKobTsu21c} to use anisotropic mesh partitions. We construct two affine simplices $\widetilde{T} \subset \mathbb{R}^d$, and two affine mappings $\Phi_{\widetilde{T}}: \widehat{T} \to \widetilde{T}$ and $\Phi_{T}: \widetilde{T} \to T$. First, we define the affine mapping $\Phi_{\widetilde{T}}: \widehat{T} \to \widetilde{T}$ as
\begin{align}
\displaystyle
\Phi_{\widetilde{T}}: \widehat{T} \ni \hat{x} \mapsto \tilde{x} := \Phi_{\widetilde{T}}(\hat{x}) := {A}_{\widetilde{T}} \hat{x} \in  \widetilde{T}, \label{aff=1}
\end{align}
where ${A}_{\widetilde{T}} \in \mathbb{R}^{d \times d}$ is an invertible matrix. We then define the affine mapping $\Phi_{T}: \widetilde{T} \to T$ as follows:
\begin{align}
\displaystyle
\Phi_{T}: \widetilde{T} \ni \tilde{x} \mapsto x := \Phi_{T}(\tilde{x}) := {A}_{T} \tilde{x} + b_{T} \in T, \label{aff=2}
\end{align}
where $b_{T} \in \mathbb{R}^d$ is a vector and ${A}_{T} \in O(d)$ denotes the rotation and mirror-imaging matrix. We define the affine mapping $\Phi: \widehat{T} \to T$ as
\begin{align*}
\displaystyle
\Phi := {\Phi}_{T} \circ {\Phi}_{\widetilde{T}}: \widehat{T} \ni \hat{x} \mapsto x := \Phi (\hat{x}) =  ({\Phi}_{T} \circ {\Phi}_{\widetilde{T}})(\hat{x}) = {A} \hat{x} + b_{T} \in T, 
\end{align*}
where ${A} := {A}_{T} {A}_{\widetilde{T}} \in \mathbb{R}^{d \times d}$.

\subsubsection*{Construct mapping $\Phi_{\widetilde{T}}: \widehat{T} \to \widetilde{T}$} \label{sec221} 
We consider the affine mapping \eqref{aff=1}. We define the matrix $ {A}_{\widetilde{T}} \in \mathbb{R}^{d \times d}$ as follows: We first define the diagonal matrix as
\begin{align}
\displaystyle
\widehat{A} :=  \diag (h_1,\ldots,h_d), \quad h_i \in \mathbb{R}_+ \quad \forall i,\label{aff=3}
\end{align}
where $\mathbb{R}_+$ denotes the set of positive real numbers.

For $d=2$, we define the regular matrix $\widetilde{A} \in \mathbb{R}^{2 \times 2}$ as
\begin{align}
\displaystyle
\widetilde{A} :=
\begin{pmatrix}
1 & s \\
0 & t \\
\end{pmatrix}, \label{aff=4}
\end{align}
with the parameters
\begin{align*}
\displaystyle
s^2 + t^2 = 1, \quad t \> 0.
\end{align*}
For the reference element $\widehat{T}$, let $\mathfrak{T}^{(2)}$ be a family of triangles.
\begin{align*}
\displaystyle
\widetilde{T} &= \Phi_{\widetilde{T}}(\widehat{T}) = {A}_{\widetilde{T}} (\widehat{T}), \quad {A}_{\widetilde{T}} := \widetilde {A} \widehat{A}
\end{align*}
with the vertices $\tilde{p}_1 := (0,0)^T$, $\tilde{p}_2 := (h_1,0)^T$, and $\tilde{p}_3 :=(h_2 s , h_2 t)^T$. Then, $h_1 = |\tilde{p}_1 - \tilde{p}_2| \> 0$ and $h_2 = |\tilde{p}_1 - \tilde{p}_3| \> 0$. 

For $d=3$, we define the regular matrices $\widetilde{A}_1, \widetilde{A}_2 \in \mathbb{R}^{3 \times 3}$ as 
\begin{align}
\displaystyle
\widetilde{A}_1 :=
\begin{pmatrix}
1 & s_1 & s_{21} \\
0 & t_1  & s_{22}\\
0 & 0  & t_2\\
\end{pmatrix}, \
\widetilde{A}_2 :=
\begin{pmatrix}
1 & - s_1 & s_{21} \\
0 & t_1  & s_{22}\\
0 & 0  & t_2\\
\end{pmatrix} \label{aff=5}
\end{align}
with the parameters
\begin{align*}
\displaystyle
\begin{cases}
s_1^2 + t_1^2 = 1, \ s_1 \> 0, \ t_1 \> 0, \ h_2 s_1 \leq h_1 / 2, \\
s_{21}^2 + s_{22}^2 + t_2^2 = 1, \ t_2 \> 0, \ h_3 s_{21} \leq h_1 / 2.
\end{cases}
\end{align*}
Therefore, we set $\widetilde{A} \in \{ \widetilde{A}_1 , \widetilde{A}_2 \}$. For the reference elements $\widehat{T}_i$, $i=1,2$, let $\mathfrak{T}_i^{(3)}$ and $i=1,2$ be a family of tetrahedra.
\begin{align*}
\displaystyle
\widetilde{T}_i &= \Phi_{\widetilde{T}_i} (\widehat{T}_i) =  {A}_{\widetilde{T}_i} (\widehat{T}_i), \quad {A}_{\widetilde{T}_i} := \widetilde {A}_i \widehat{A}, \quad i=1,2,
\end{align*}
with the vertices
\begin{align*}
\displaystyle
&\tilde{p}_1 := (0,0,0)^T, \ \tilde{p}_2 := (h_1,0,0)^T, \ \tilde{p}_4 := (h_3 s_{21}, h_3 s_{22}, h_3 t_2)^T, \\
&\begin{cases}
\tilde{p}_3 := (h_2 s_1 , h_2 t_1 , 0)^T \quad \text{for case (\roman{sone})}, \\
\tilde{p}_3 := (h_1 - h_2 s_1, h_2 t_1,0)^T \quad \text{for case (\roman{stwo})}.
\end{cases}
\end{align*}
Subsequently, $h_1 = |\tilde{p}_1 - \tilde{p}_2| \> 0$, $h_3 = |\tilde{p}_1 - \tilde{p}_4| \> 0$, and
\begin{align*}
\displaystyle
h_2 =
\begin{cases}
|\tilde{p}_1 - \tilde{p}_3| \> 0  \quad \text{for case (\roman{sone})}, \\
|\tilde{p}_2 - \tilde{p}_3| \> 0  \quad \text{for case (\roman{stwo})}.
\end{cases}
\end{align*}

\subsubsection*{Construct mapping $\Phi_{T}: \widetilde{T} \to T$}  \label{sec322}
We determine the affine mapping \eqref{aff=2} as follows: Let ${T} \in \mathbb{T}_h$ have vertices ${p}_i$ ($i=1,\ldots,d+1$). Let $b_{T} \in \mathbb{R}^d$ be the vector and ${A}_{T} \in O(d)$ be the rotation and mirror imaging matrix such that
\begin{align*}
\displaystyle
p_{i} = \Phi_T (\tilde{p}_i) = {A}_{T} \tilde{p}_i + b_T, \quad i \in \{1, \ldots,d+1 \},
\end{align*}
where vertices $p_{i}$ ($i=1,\ldots,d+1$) satisfy the following conditions:

\begin{Cond}[Case in which $d=2$] \label{cond1}
Let ${T} \in \mathbb{T}_h$ have vertices ${p}_i$ ($i=1,\ldots,3$). We assume that $\overline{{p}_2 {p}_3}$ is the longest edge of ${T}$; i.e. $ h_{{T}} := |{p}_2 - {p}_ 3|$. We set $h_1 = |{p}_1 - {p}_2|$ and $h_2 = |{p}_1 - {p}_3|$. Then, we assume that $h_2 \leq h_1$. Note that ${h_1 \approx h_T}$. 
\end{Cond}

\begin{Cond}[Case in which $d=3$] \label{cond2}
Let ${T} \in \mathbb{T}_h$ have vertices ${p}_i$ ($i=1,\ldots,4$). Let ${L}_i$ ($1 \leq i \leq 6$) be the edges of ${T}$. We denote by ${L}_{\min}$  the edge of ${T}$ with the minimum length; i.e. $|{L}_{\min}| = \min_{1 \leq i \leq 6} |{L}_i|$. We set $h_2 := |{L}_{\min}|$ and assume that 
\begin{align*}
\displaystyle
&\text{the endpoints of ${L}_{\min}$ are either $\{ {p}_1 , {p}_3\}$ or $\{ {p}_2 , {p}_3\}$}.
\end{align*}
Among the four edges that share an endpoint with ${L}_{\min}$, we consider the longest edge ${L}^{({\min})}_{\max}$. Let ${p}_1$ and ${p}_2$ be the endpoints of edge ${L}^{({\min})}_{\max}$. Thus, we have
\begin{align*}
\displaystyle
h_1 = |{L}^{(\min)}_{\max}| = |{p}_1 - {p}_2|.
\end{align*}
We consider cutting $\mathbb{R}^3$ with the plane that contains the midpoint of the edge ${L}^{(\min)}_{\max}$ and is perpendicular to the vector ${p}_1 - {p}_2$. Thus, there are two cases. 
\begin{description}
  \item[(Type \roman{sone})] ${p}_3$ and ${p}_4$  belong to the same half-space;
  \item[(Type \roman{stwo})] ${p}_3$ and ${p}_4$  belong to different half-spaces.
\end{description}
In each case, we set
\begin{description}
  \item[(Type \roman{sone})] ${p}_1$ and ${p}_3$ as the endpoints of ${L}_{\min}$, that is, $h_2 =  |{p}_1 - {p}_3| $;
  \item[(Type \roman{stwo})] ${p}_2$ and ${p}_3$ as the endpoints of ${L}_{\min}$, that is, $h_2 =  |{p}_2 - {p}_3| $.
\end{description}
Finally, we set $h_3 = |{p}_1 - {p}_4|$. We implicitly assume that ${p}_1$ and ${p}_4$ belong to the same half-space. In addition, note that ${h_1 \approx h_T}$. 
\end{Cond}

\subsection{Additional notations and assumptions} \label{addinot}
\begin{figure}[htbp]
  \includegraphics[keepaspectratio, scale=0.45]{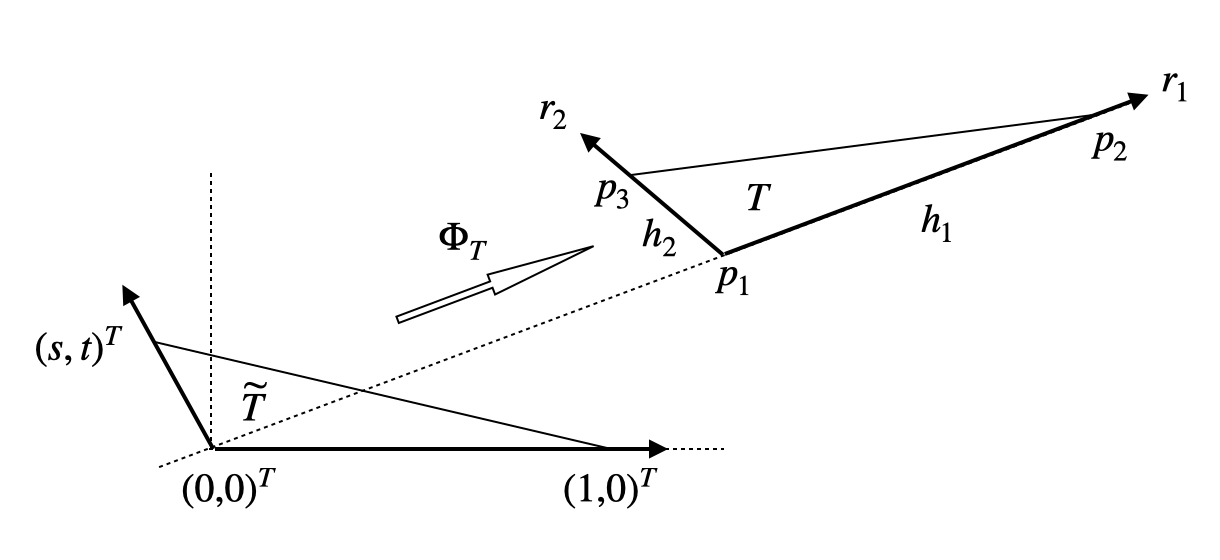}
\caption{Affine mapping $\Phi_{T}$ and vectors $r_i$, $i=1,2$}
\label{affine_2d}
\end{figure}

\begin{figure}[htbp]
  \begin{minipage}[b]{0.4\linewidth}
    \centering
    \includegraphics[keepaspectratio, scale=0.45]{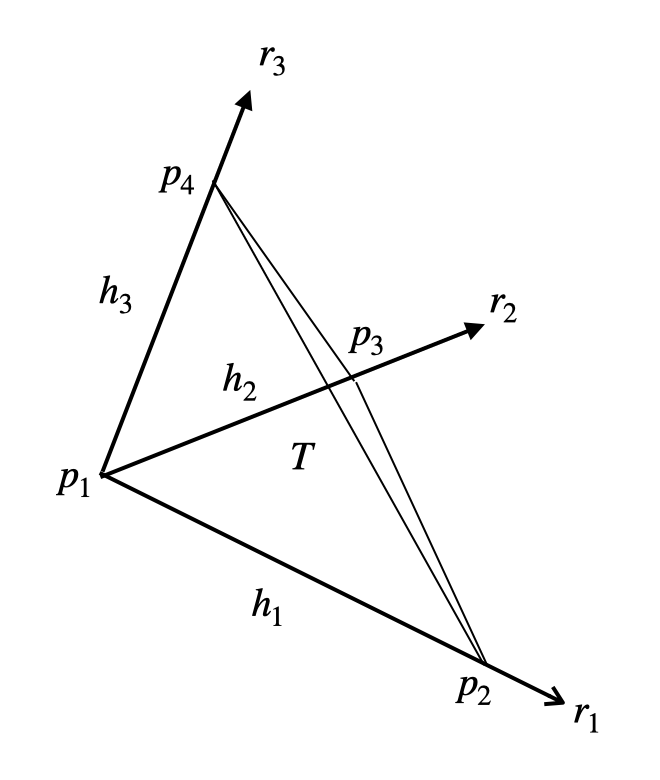}
    \caption{(Type \roman{sone}) Vectors $r_i$, $i=1,2,3$}
     \label{affine_3d_1}
  \end{minipage}
  \begin{minipage}[b]{0.4\linewidth}
    \centering
    \includegraphics[keepaspectratio, scale=0.45]{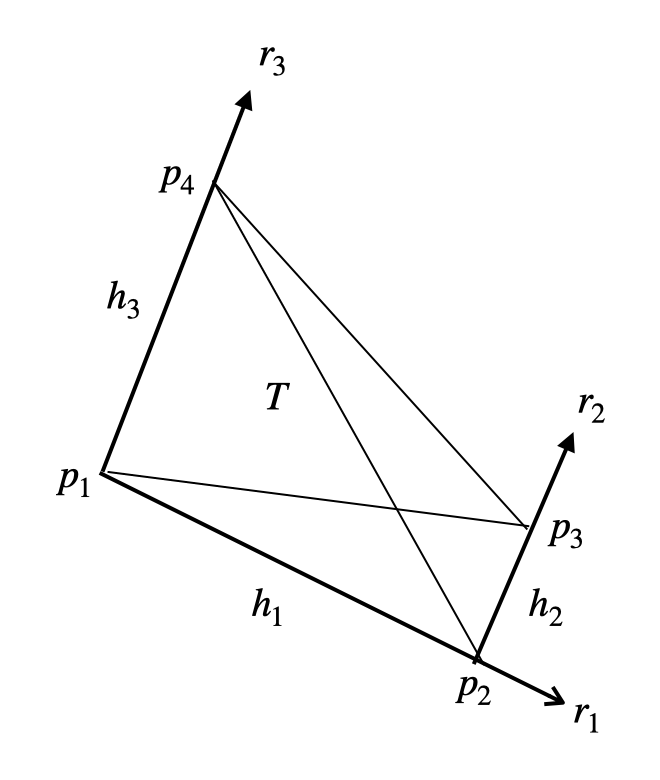}
    \caption{(Type \roman{stwo}) Vectors $r_i$, $i=1,2,3$}
     \label{affine_3d_2}
  \end{minipage}
\end{figure}

We define the vectors ${r}_n \in \mathbb{R}^d$, $n=1,\ldots,d$ as follows: If $d=2$,
\begin{align*}
\displaystyle
{r}_1 := \frac{p_2 - p_1}{|p_2 - p_1|}, \quad {r}_2 := \frac{p_3 - p_1}{|p_3 - p_1|},
\end{align*}
see Fig. \ref{affine_2d}, and if $d=3$,
\begin{align*}
\displaystyle
&{r}_1 := \frac{p_2 - p_1}{|p_2 - p_1|}, \quad {r}_3 := \frac{p_4 - p_1}{|p_4 - p_1|}, \quad
\begin{cases}
\displaystyle
{r}_2 := \frac{p_3 - p_1}{|p_3 - p_1|}, \quad \text{for (Type \roman{sone})}, \\
\displaystyle
{r}_2 := \frac{p_3 - p_2}{|p_3 - p_2|} \quad \text{for (Type \roman{stwo})},
\end{cases}
\end{align*}
see Fig \ref{affine_3d_1} for (Type \roman{sone}) and Fig \ref{affine_3d_2} for (Type \roman{stwo}).

For a sufficiently smooth function $\varphi$ and a vector function $v := (v_{1},\ldots,v_{d})^T$, we define the directional derivative for $i \in \{ 1, \ldots,d \}$:
\begin{align*}
\displaystyle
\frac{\partial \varphi}{\partial {r_i}} &:= ( {r}_i \cdot  \nabla_{x} ) \varphi = \sum_{i_0=1}^d ({r}_i)_{i_0} \frac{\partial \varphi}{\partial x_{i_0}^{}}, \\
\frac{\partial v}{\partial r_i} &:= \left(\frac{\partial v_{1}}{\partial r_i}, \ldots, \frac{\partial v_{d}}{\partial r_i} \right)^T 
= ( ({r}_i  \cdot \nabla_{x}) v_{1}, \ldots, ({r}_i  \cdot \nabla_{x} ) v_{d} )^T.
\end{align*}
For a multiindex $\beta = (\beta_1,\ldots,\beta_d) \in \mathbb{N}_0^d$, we use the notation
\begin{align*}
\displaystyle
\partial^{\beta} \varphi := \frac{\partial^{|\beta|} \varphi}{\partial x_1^{\beta_1} \ldots \partial x_d^{\beta_d}}, \quad \partial^{\beta}_{r} \varphi := \frac{\partial^{|\beta|} \varphi}{\partial r_1^{\beta_1} \ldots \partial r_d^{\beta_d}}, \quad h^{\beta} :=  h_{1}^{\beta_1} \cdots h_{d}^{\beta_d}.
\end{align*}
Note that $\partial^{\beta} \varphi \neq  \partial^{\beta}_{r} \varphi$.

We proposed a geometric parameter $H_{T}$ in a prior paper \cite{IshKobTsu21a}.
 \begin{defi} \label{defi1}
 The parameter $H_{{T}}$ is defined as
\begin{align*}
\displaystyle
H_{{T}} := \frac{\prod_{i=1}^d h_i}{|{T}|_d} h_{{T}}.
\end{align*}
\end{defi}
Here, we introduce the geometric condition proposed in \cite{IshKobTsu21a}, which is equivalent to the maximum angle condition \cite{IshKobSuzTsu21d}.

\begin{assume} \label{neogeo=assume}
A family of meshes $\{ \mathbb{T}_h\}$ has a semi-regular property if there exists $\gamma_0 \> 0$ such that
\begin{align}
\displaystyle
\frac{H_{T}}{h_{T}} \leq \gamma_0 \quad \forall \mathbb{T}_h \in \{ \mathbb{T}_h \}, \quad \forall T \in \mathbb{T}_h. \label{NewGeo}
\end{align}
\end{assume}

{
\begin{rem}
In our anisotropic interpolation errors, the quantity $ \| \widetilde{{A}} \|_2 \| \widetilde{{A}}^{-1} \|_2 $ is estimated as follows:
\begin{align*}
\displaystyle
 \| \widetilde{{A}} \|_2 \| \widetilde{{A}}^{-1} \|_2 \leq c \frac{H_{T}}{h_{T}},
\end{align*}
see \cite[Lemma 2]{IshKobTsu21c} for the proof. Therefore, by imposing Assumption \ref{neogeo=assume}, the quantity $ \| \widetilde{{A}} \|_2 \| \widetilde{{A}}^{-1} \|_2 $ is bounded. 	
\end{rem}
}

\subsection{$L^2$-orthogonal projection}
For $T_j \in \mathbb{T}_h$, $j \in \{ 1, \ldots , Ne \}$, let $\Pi_{T_j}^0 : L^2(T_j) \to \mathbb{P}^0(T_j)$ be the $L^2$-orthogonal projection defined as
\begin{align*}
\displaystyle
\Pi_{T_j}^0 \varphi := \frac{1}{|T_j|_d} \int_{{T_j}} \varphi dx \quad \forall \varphi \in L^2(T_j).
\end{align*}
The following theorem provides an anisotropic error estimate for the projection $\Pi_{T_j}^0$. 
\begin{thr} \label{thr1}
For any $\hat{\varphi} \in H^{1}(\widehat{T})$ with ${\varphi} := \hat{\varphi} \circ {\Phi}^{-1}$,
\begin{align}
\displaystyle
\| \Pi_{T_j}^0 \varphi - \varphi \|_{L^2(T_j)} \leq c \sum_{i=1}^d h_i \left\| \frac{\partial \varphi}{\partial r_i} \right\|_{L^{2}(T_j)}. \label{L2ortho}
\end{align}
\end{thr}

\begin{pf*}
The proof can be found in \cite[Theorem 2]{Ish24} and \cite[Theorem 2]{Ish24b}.
\qed
\end{pf*}

We also define the global interpolation $\Pi_h^0$ to space $P_{dc,h}^{0}$ as
\begin{align*}
\displaystyle
(\Pi_h^0 \varphi)|_{T_j} := \Pi_{{T_j}}^0 (\varphi|_{T_j}) \ \forall T_j \in \mathbb{T}_h, \ j \in \{ 1, \ldots ,Ne\}, \ \forall \varphi \in L^2(\Omega).
\end{align*}

\subsection{CR finite element interpolation operator}
For $T_j \in \mathbb{T}_h$, $j \in \{ 1, \ldots , Ne \}$, let the points $\{ P_{T_j,1}, \ldots, P_{T_j,d+1} \}$ be the vertices of the simplex $T_j \in \mathbb{T}_h$. Let $F_{T_j,i}$ be the face of $T_j$ opposite $P_{T_j,i}$ for $i \in \{ 1, \ldots , d+1\}$. The CR interpolation operator $I_{T_j}^{CR} : H^{1}(T_j) \to \mathbb{P}^1(T_j)$ is defined as, for any $\varphi \in H^{1}(T_j)$, 
\begin{align*}
\displaystyle
I_{T_j}^{CR}: H^{1}(T_j) \ni \varphi \mapsto I_{T_j}^{CR} \varphi := \sum_{i=1}^{d+1} \left(  \frac{1}{| {F}_{T_j,i} |_{d-1}} \int_{{F}_{T_j,i}} {\varphi} d{s} \right) \theta_{T_j,i}^{CR} \in \mathbb{P}^1(T_j). 
\end{align*}
where the basis of the CR finite element is defined in \eqref{CR2}. We then present the estimates of the anisotropic CR interpolation error. 

\begin{thr} \label{DGCE=thr3}
For $j \in \{ 1, \ldots , Ne\}$,
\begin{align}
\displaystyle
|I_{T_j}^{CR} \varphi - \varphi |_{H^{1}({T}_j)} &\leq c \sum_{i=1}^d {h}_i \left\| \frac{\partial }{\partial r_i} \nabla \varphi \right \|_{L^2(T_j)^d} \quad \forall {\varphi} \in H^{2}({T}_j), \label{CR4} \\
\|I_{T_j}^{CR} \varphi - \varphi \|_{L^2(T_j)}
&\leq c  \sum_{|\varepsilon| = 2} h^{\varepsilon} \left\| \partial_{r}^{\varepsilon} \varphi  \right\|_{L^{2}(T)} \quad \forall {\varphi} \in H^{2}({T}_j). \label{L2=CR4}
\end{align}
\end{thr}

\begin{pf*}
The proof of \eqref{CR4} can be found in \cite[Theorem 3]{Ish24} and \cite[Theorem 3]{Ish24b}. The proof of \eqref{L2=CR4} can be found in Appendix.
\qed
\end{pf*}

We define the global interpolation operator ${I}_{h}^{CR}: H^{1}(\Omega) \to V_{dc,h}^{CR}$ as
\begin{align}
\displaystyle
({I}_{h}^{CR} \varphi )|_{T_j} = {I}_{T_j}^{CR} (\varphi |_{T_j}), \quad j \in \{ 1, \ldots , Ne\}, \quad \forall \varphi \in H^{1}(\Omega), \label{CR8}
\end{align}	
where the space $V_{dc,h}^{CR}$ is defined in \eqref{CR6}.

\subsection{RT finite element interpolation operator} \label{RTsp}
For $T_j \in \mathbb{T}_h$, $j \in \{ 1, \ldots , Ne \}$, we define the local RT polynomial space as follows:
\begin{align}
\displaystyle
\mathbb{RT}^0(T_j) := \mathbb{P}^0(T_j)^d + x \mathbb{P}^0(T_j), \quad x \in \mathbb{R}^d. \label{RT1}
\end{align}
Let $\mathcal{I}_{T_j}^{RT}: H^{1}(T_j)^d \to \mathbb{RT}^0(T_j)$ be the RT interpolation operator such that for any $v \in H^{1}(T_j)^d$,
\begin{align}
\displaystyle
\mathcal{I}_{T_j}^{RT}: H^{1}(T_j)^d \ni v \mapsto \mathcal{I}_{T_j}^{RT} v := \sum_{i=1}^{d+1} \left(  \int_{{F}_{T_j,i}} {v} \cdot n_{T_j,i} d{s} \right) \theta_{T_j,i}^{RT} \in \mathbb{RT}^0(T_j), \label{RT4}
\end{align}
where $\theta_{T_j,i}^{RT}$ is the local shape function (e.g. \cite{ErnGue21a}) and $n_{T_j,i}$ is a fixed unit normal to ${F}_{T_j,i}$. 

The following two theorems are divided into elements of (Type \roman{sone}) or (Type \roman{stwo}) in Section \ref{two=step} when $d=3$.

\begin{thr} \label{DGRT=thr3}
Let $T_j$ be an element with Conditions \ref{cond1} or \ref{cond2} satisfying (Type \roman{sone}) in Section \ref{two=step} when $d=3$. For any ${v} = ({v}_1,\ldots,{v}_d)^T \in H^1(T_j)^d$ and $j \in \{ 1, \ldots ,Ne\}$,
\begin{align}
\displaystyle
\| \mathcal{I}_{T_j}^{RT} v - v \|_{L^2(T_j)^d} 
&\leq  c \left( \frac{H_{T_j}}{h_{T_j}} \sum_{i=1}^d h_i \left \|  \frac{\partial v}{\partial r_i} \right \|_{L^2(T_j)^d} +  h_{T_j} \| \div {v} \|_{L^{2}({T}_j)} \right). \label{RT5}
\end{align}
\end{thr}

\begin{pf*}
The proof is provided in \cite[Theorem 2]{Ish21}.
\qed
\end{pf*}

\begin{thr} \label{DGRT=thr4}
Let $d=3$. Let $T_j$ be an element with Condition \ref{cond2} that satisfies (Type \roman{stwo}) in Section \ref{two=step}. For ${v} = ({v}_1,v_2,{v}_3)^T \in H^1(T_j)^3$ and $j \in \{ 1, \ldots ,Ne\}$,
\begin{align}
\displaystyle
&\| \mathcal{I}_{T_j}^{RT} v - v \|_{L^2(T_j)^3} 
\leq c \frac{H_{T_j}}{h_{T_j}} \Biggl(  h_{T_j} |v|_{H^1(T_j)^3} \Biggr). \label{RT6}
\end{align}
\end{thr}

\begin{pf*}
The proof is provided in \cite[Theorem 3]{Ish21}.
\qed
\end{pf*}

A discontinuous RT finite element space is defined as follows:
\begin{align}
\displaystyle
V_{dc,h}^{RT} &:= \{ v_h \in L^1(\Omega)^d : \  v_h|_{T_j} \in \mathbb{RT}^0({T}) \quad \forall T_j \in \mathbb{T}_h, \  j \in \{ 1, \ldots , Ne\} \}. \label{RT7}
\end{align}
We also define a global interpolation operator $\mathcal{I}_{h}^{RT}: H^{1}(\Omega)^d  \to V_{dc,h}^{RT}$ as
\begin{align}
\displaystyle
(\mathcal{I}_{h}^{RT} v )|_{T_j} = \mathcal{I}_{T_j}^{RT} (v |_{T_j}), \quad j \in \{ 1, \ldots , Ne\}, \quad \forall v \in H^{1}(\Omega)^d.  \label{RT8}
\end{align}

\section{Discrete Poincar\'e inequality} \label{disPoi}
In this section, we present the discrete Poincar\'e inequality for anisotropic meshes. Here, we present a proof that uses a dual problem. Therefore, we impose convexity on $\Omega$ to establish the aforementioned inequality.

Now, we introduce the Jensen-type inequality (see \cite[Exercise 12.1]{ErnGue21a}). Let $r,s$ be two nonnegative real numbers and $\{ x_i \}_{i \in I}$ be a finite sequence of nonnegative numbers. Thus, it holds that
\begin{align}
\displaystyle
\begin{cases}
\left( \sum_{i \in I} x_i^s \right)^{\frac{1}{s}} \leq \left( \sum_{i \in I} x_i^r \right)^{\frac{1}{r}} \quad \text{if $r \leq s$},\\
\left( \sum_{i \in I} x_i^s \right)^{\frac{1}{s}} \leq \card(I)^{\frac{r-s}{rs}} \left( \sum_{i \in I} x_i^r \right)^{\frac{1}{r}} \quad \text{if $r \> s$}.
\end{cases} \label{jensen}
\end{align}

The following relation plays an important role in the discontinuous Galerkin finite element analysis on anisotropic meshes.

\begin{lem} \label{rel=CRRT}
For any $w \in H^1(\Omega)^d$, $\psi_h \in P_{dc,h}^{1}$ and $\lambda_h \in V_{dc,h,0}(\mathcal{F}_h)$,
\begin{align}
\displaystyle
&\int_{\Omega} \left( \mathcal{I}_h^{RT} w \cdot \nabla_h \psi_{h} + \div \mathcal{I}_h^{RT} w  \psi_{h} \right) dx \notag\\
&\quad  = \sum_{T \in \mathbb{T}_h} \sum_{F \in \mathcal{F}_T} n_{T,F} \int_{F} ( w \cdot n_F) \Pi_F^0 (\psi_{h} - \lambda_h ) ds,  \label{hwop=1}
\end{align}
where $n_{T,F} := n_T \cdot n_F$.
\end{lem}

\begin{pf*}
For any $w \in H^1(\Omega)^d$ and $\psi_h \in P_{dc,h}^{1}$, using Green's formula, we derive
\begin{align*}
\displaystyle
&\int_{\Omega} \left( \mathcal{I}_h^{RT} w \cdot \nabla_h \psi_{h} + \div \mathcal{I}_h^{RT} w  \psi_{h} \right) dx = \sum_{T \in \mathbb{T}_h} \sum_{F \in \mathcal{F}_T} n_{T,F} \int_{F} (\mathcal{I}_h^{RT} w \cdot n_F) \psi_{h} ds.
\end{align*}
Because $\lambda_h \in V_{dc,h,0}(\mathcal{F}_h)$ is single-valued on $F \in \mathcal{F}_h$ and $\mathcal{I}_h^{RT} w \cdot n_F \in \mathbb{P}^{0}(F)$, for any $F \in \mathcal{F}_h$,
\begin{align*}
\displaystyle
\sum_{T \in \mathbb{T}_h} \sum_{F \in \mathcal{F}_T} n_{T,F} \int_{F} (\mathcal{I}_h^{RT} w \cdot n_F) \lambda_{h} ds = 0.
\end{align*}
The properties of the projection $\Pi_F^0$ and $\mathcal{I}_h^{RT}$ yield
\begin{align*}
\displaystyle
&\sum_{T \in \mathbb{T}_h} \sum_{F \in \mathcal{F}_T} n_{T,F} \int_{F} (\mathcal{I}_h^{RT} w \cdot n_F) \psi_{h} ds \\
&\quad = \sum_{T \in \mathbb{T}_h} \sum_{F \in \mathcal{F}_T} n_{T,F} \int_{F} (\mathcal{I}_h^{RT} w \cdot n_F) (\psi_{h} - \lambda_h ) ds \\
&\quad = \sum_{T \in \mathbb{T}_h} \sum_{F \in \mathcal{F}_T} n_{T,F} \int_{F} (\mathcal{I}_h^{RT} w \cdot n_F) \Pi_F^0 (\psi_{h} - \lambda_h ) ds \\
&\quad = \sum_{T \in \mathbb{T}_h} \sum_{F \in \mathcal{F}_T} n_{T,F} \int_{F} ( w \cdot n_F) \Pi_F^0 (\psi_{h} - \lambda_h ) ds,
\end{align*}
which leads to the target equality \eqref{hwop=1}.
\qed
\end{pf*}

The right-hand terms in \eqref{hwop=1} are estimated as follows:

\begin{lem} \label{est=CRRT}
Let $\beta$ be nonnegative real numbers. For any $w \in H^1(\Omega)^d$, $\psi_h \in P_{dc,h}^{1}$ and $\lambda_h \in V_{dc,h,0}(\mathcal{F}_h)$,
\begin{align}
\displaystyle
&\left| \sum_{T \in \mathbb{T}_h} \sum_{F \in \mathcal{F}_T} n_{T,F} \int_{F} ( w \cdot n_F) \Pi_F^0 (\psi_{h} - \lambda_h ) ds \right| \notag\\
&\quad \leq c \left(  h^{\beta} \| w \|_{L^2(\Omega)^d} + h^{\beta + \frac{1}{2}} \| w \|_{L^2(\Omega)^d}^{\frac{1}{2}} | w |_{H^1(\Omega)^d}^{\frac{1}{2}} \right) \notag \\
&\quad \quad \times  \left(  \sum_{T \in \mathbb{T}_h} \sum_{F \in \mathcal{F}_T} h^{- 2 \beta} \ell_{T,F}^{- 1} \left \|  \Pi_F^0 (\psi_{h} - \lambda_h ) \right\|_{L^2(F)}^2 \right)^{\frac{1}{2}}. \label{hwop=2}
\end{align}
\end{lem}

\begin{pf*}
Using the H\"older and Cauchy--Schwarz inequalities and the trace inequality \eqref{trace=vec}, we obtain
\begin{align*}
\displaystyle
&\left| \sum_{T \in \mathbb{T}_h} \sum_{F \in \mathcal{F}_T} n_{T,F} \int_{F} ( w \cdot n_F) \Pi_F^0 (\psi_{h} - \lambda_h ) ds \right| \\
&\quad \leq c \sum_{T \in \mathbb{T}_h} \sum_{F \in \mathcal{F}_T} \ell_{T,F}^{- \frac{1}{2}} \left( \| w \|_{L^2(T)^d} + h_{T}^{\frac{1}{2}}  \| w \|_{L^2(T)^d}^{\frac{1}{2}} | w |_{H^1(T)^d}^{\frac{1}{2}} \right) \left \|  \Pi_F^0 (\psi_{h} - \lambda_h ) \right\|_{L^2(F)} \\
&\quad \leq c \left(  \sum_{T \in \mathbb{T}_h} h^{2 \beta} \left( \| w \|_{L^2(T)^d} + h_{T}^{\frac{1}{2}}  \| w \|_{L^2(T)^d}^{\frac{1}{2}} | w |_{H^1(T)^d}^{\frac{1}{2}} \right)^2 \right)^{\frac{1}{2}} \\
&\quad \quad \times  \left(  \sum_{T \in \mathbb{T}_h} \sum_{F \in \mathcal{F}_T} h^{- 2 \beta} \ell_{T,F}^{- 1} \left \|  \Pi_F^0 (\psi_{h} - \lambda_h ) \right\|_{L^2(F)}^2 \right)^{\frac{1}{2}},
\end{align*}
which leads to the inequality \eqref{hwop=2} together with the Jensen-type inequality \eqref{jensen}.
\qed
\end{pf*}

The following lemma provides the discrete Poincar\'e inequality. For simplicity, we assume that $\Omega$ is convex.

\begin{lem}[Discrete Poincar\'e inequality] \label{DisPoi=lem}
Assume that $\Omega$ is a convex. Let $\{ \mathbb{T}_h\}$ be a family of meshes with semi-regular properties (Assumption \ref{neogeo=assume}) and $h \leq 1$. Then, there exists a positive constant $C_{dc}^{P}$ independent of $h$ {but dependent on the maximum angle}, such that
\begin{align}
\displaystyle
 \| \varphi_h \|_{L^2(\Omega)} \leq C_{dc}^{P} | \varphi_h^H |_{hwop(0)} \quad \forall \varphi_h^H := (\varphi_h,\mu_h) \in V_{dc,h,0}^{H}. \label{hwop=3}
\end{align}
\end{lem}

\begin{pf*}
Let $\varphi_h^H := (\varphi_h,\mu_h) \in V_{dc,h,0}^{H}$. Consider the following problem: We determine $z \in H^{2}(\Omega) \cap H_0^1(\Omega)$ such that:
\begin{align*}
\displaystyle
- \varDelta z = \varphi_{h} \quad \text{in $\Omega$}, \quad z = 0 \quad \text{on $\partial \Omega$}.
\end{align*}
Then, we obtain a priori estimates $|z|_{H^1(\Omega)} \leq C_P \| \varphi_{h} \|_{L^2(\Omega)}$ and $ |z|_{H^2(\Omega)} \leq \| \varphi_h \|_{L^2(\Omega)}$, where $C_P$ is the Poincar\'e constant.

We provide the following {equality} for the analysis:
\begin{align*}
\displaystyle
 \int_{\Omega} ( \div \mathcal{I}_{h}^{RT}  (\nabla z) ) \varphi_h {dx}
 &= \sum_{T \in \mathbb{T}_h} \int_{\partial T} \mathcal{I}_{h}^{RT} (\nabla z) \cdot n_T \varphi_h ds - \int_{\Omega} \mathcal{I}_{h}^{RT} (\nabla z) \cdot \nabla_h \varphi_h dx \\
 &=  \int_{\Omega} ( \nabla z -  \mathcal{I}_{h}^{RT}  (\nabla z) ) \cdot \nabla_h \varphi_h dx  - \int_{\Omega}  \nabla z \cdot \nabla_h \varphi_h dx \\
 &\quad + \sum_{T \in \mathbb{T}_h} \int_{\partial T} \mathcal{I}_{h}^{RT}  (\nabla z) \cdot n_T \varphi_h ds \\
  &=  \int_{\Omega} ( \nabla z -  \mathcal{I}_{h}^{RT}  (\nabla z) ) \cdot \nabla_h \varphi_h dx  - \int_{\Omega}  \nabla z \cdot \nabla_h \varphi_h dx \\
  &\quad +  \sum_{T \in \mathbb{T}_h} \sum_{F \in \mathcal{F}_T} n_{T,F} \int_{F} ( \nabla z \cdot n_F) \Pi_F^0 (\varphi_{h} - \mu_h ) ds,
 \end{align*}
where we apply an analogous argument to Lemma \ref{rel=CRRT} for the last equality. { Equality } yields
\begin{align*}
\displaystyle
&\| \varphi_h \|^2_{L^2(\Omega)} = \int_{\Omega} \varphi_h^2 dx = \int_{\Omega} - \varDelta z \varphi_h dx = - \int_{\Omega} \div (\nabla z)  \varphi_h dx \\
&\quad =  \int_{\Omega} ( \Pi_h^{0}  \div (\nabla z) - \div (\nabla z) ) \varphi_h dx - \int_{\Omega} ( \Pi_h^{0}  \div (\nabla z) ) \varphi_h dx\\
&\quad =  \int_{\Omega} ( \Pi_h^{0}  \div (\nabla z) - \div (\nabla z) ) ( \varphi_h - \Pi_h^{0} \varphi_h) dx - \int_{\Omega} ( \div \mathcal{I}_{h}^{RT}  (\nabla z) ) \varphi_h dx\\
&\quad =  - \int_{\Omega} \div (\nabla z) \left( \psi_h - \Pi_h^{0} \varphi_h \right) dx \\
&\quad \quad  - \int_{\Omega} ( \nabla z -  \mathcal{I}_{h}^{RT}  (\nabla z) ) \cdot \nabla_h \varphi_h dx  + \int_{\Omega}  \nabla z \cdot \nabla_h \varphi_h dx \\
&\quad \quad -  \sum_{T \in \mathbb{T}_h} \sum_{F \in \mathcal{F}_T} n_{T,F} \int_{F} ( \nabla z \cdot n_F) \Pi_F^0 (\varphi_{h} - \mu_h ) ds.
\end{align*}
 Using H\"older's inequality, the error estimates \eqref{L2ortho}, \eqref{RT5}, and \eqref{RT6}, and Lemma \ref{est=CRRT} with $\beta = 0$, we have
\begin{align*}
\displaystyle
\| \varphi_h \|^2_{L^2(\Omega)}
&\leq \| \varDelta z \|_{L^2(\Omega)} \| \varphi_h - \Pi_h^{0} \varphi_h \|_{L^2(\Omega)} + \|  \nabla z -  \mathcal{I}_{h}^{RT}  (\nabla z) \|_{L^2(\Omega)} |\varphi_h|_{H^1(\mathbb{T}_h)} \\
&\quad + |z|_{H^1(\Omega)} |\varphi_h|_{H^1(\mathbb{T}_h)} + c \| \nabla z \|_{H^1(\Omega)^d} |\varphi_h^H|_{hwop(0)} \\
&\leq c( h + 1 ) \| \varphi_h \|_{L^2(\Omega)} |\varphi_h^H|_{hwop(0)},
\end{align*}
which leads to the target inequality if $h \leq 1$.
\qed
\end{pf*}

\section{Stability and error estimates for the HWOPSIP method}

\subsection{Stability}

\begin{thr}[Stability]
Assume that $\Omega$ is convex. Let $\{ \mathbb{T}_h\}$ be a family of meshes with semiregular properties (Assumption \ref{neogeo=assume}) and $h \leq 1$. Let $u_h^H := (u_h^{hwop}, \lambda_h^{hwop}) \in V_{dc,h,0}^{H}$ be a discrete solution of \eqref{hwop=4}. Then,
\begin{align}
\displaystyle
 | u_h^H |_{hwop(1)} \leq c \| f \|_{L^2(\Omega)}. \label{hwop=7}
\end{align}
\end{thr}

\begin{pf*}
The H\"older and discrete Poincar\'e \eqref{hwop=3} inequalities yield
\begin{align*}
\displaystyle
|\ell_h(\varphi_{h}^H)| \leq C_{dc}^{P} \| f \|_{L^2(\Omega)} | \varphi_{h}^H |_{hwop(0)} \leq c  \| f \|_{L^2(\Omega)} | \varphi_{h}^H |_{hwop(1)}  \quad \forall \varphi_{h}^H \in V_{dc,h,0}^{H},
\end{align*}
if $\Omega$ is convex, and $h \leq 1$. Using \eqref{hwop=6}, we obtain the stability estimate \eqref{hwop=7}:
\qed
\end{pf*}

\begin{rem} \label{poin=rem1}
To prove the stability estimates of the schemes, we consider the case in which $\Omega$ is not convex. In \cite{Bre03}, discrete Poincar\'e inequalities for piecewise $H^1$ functions are proposed. However, the inverse, trace inequalities, and the local quasi-uniformity for meshes under the shape-regularity condition are used for the proof. Therefore, a careful consideration of the results in \cite{Bre03} may be necessary to eliminate the assumption that $\Omega$ is convex.
\end{rem}

\subsection{Energy norm error estimate}
The starting point for error analysis is the second Strang lemma; see \cite[Lemma 2.25]{ErnGue04}.

\begin{lem} \label{second=starang}
We assume that $\Omega$ is convex. Let $u \in H_0^1(\Omega)$ be the solution to \eqref{poisson_weak} and $u^H := (u,u|_{\mathcal{F}_h}) \in X_0^H$. Let $u_h^H := (u_h^{hwop}, \lambda_h^{hwop}) \ \in V_{dc,h,0}^{H}$ be the solution to Equation \eqref{hwop=4}. Then, it holds that
\begin{align}
\displaystyle
|u^H - u_h^{H}|_{hwop(1)}
&\leq \inf_{v_h^H = (v_h,\lambda_h) \in V_{dc,h,0}^{H}} |u^H - v_h^H|_{hwop(1)} + E_h(u^H), \label{hwop=8}
\end{align}
where
\begin{align}
\displaystyle
E_h(u^H) := \sup_{w_h^H = (w_h,\mu_h) \in V_{dc,h,0}^{H}} \frac{|a_h^{hwop}(u^H,w_h^H) - \ell_h(w_h^H) |}{|w_h^H|_{hwop(1)}}.  \label{hwop=9}
\end{align}
\end{lem}

\begin{pf*}
Using \eqref{hwop=5}, we have
\begin{align*}
\displaystyle
|v_h^H - u_h^{H}|_{hwop(1)}^2
&= a_h^{hwop}(v_h^H - u_h^{H},v_h^H - u_h^{H }) \\
&\hspace{-3.0cm} = a_h^{hwop}(v_h^H - u^H,v_h^H - u_h^{H}) + a_h^{hwop}(u^H,v_h^H - u_h^{H}) - \ell_h(v_h^H - u_h^{H}) \\
&\hspace{-3.0cm} \leq c |u^H - v_h^H|_{hwop(1)}  |v_h^H - u_h^{H}|_{hwop(1)} + | a_h^{hwop}(u^H,v_h^H - u_h^{H}) - \ell_h(v_h^H - u_h^{H}) |,
\end{align*}
which leads to
\begin{align*}
\displaystyle
|v_h^H - u_h^{H}|_{hwop(1)}
&\leq c |u^H - v_h^H|_{hwop(1)} + \frac{| a_h^{hwop}(u^H,v_h^H - u_h^{H}) - \ell_h(v_h^H - u_h^{H}) |}{|v_h^H - u_h^{H}|_{hwop(1)}} \\
&\leq  c |u^H - v_h^H|_{hwop(1)} + E_h(u^H).
\end{align*}
Then,
\begin{align*}
\displaystyle
|u - u_h^{H}|_{hwop(1)}
&\leq |u - v_h^H|_{hwop(1)} + |v_h^H - u_h^{H}|_{hwop(1)} \\
&\leq c  |u - v_h^H|_{hwop(1)}  + E_h(u^H).
\end{align*}
Hence, the target inequality \eqref{hwop=8} holds.
\qed
\end{pf*}

For any $T \in \mathbb{T}_h$ and $F \in \mathcal{F}_T$, we define the local discrete space as
\begin{align*}
\displaystyle
V_{T,F}^{H} := \mathbb{P}^1(T) \times  \mathbb{P}^0(F).
\end{align*}
We define the local operator $I_{T,F}^H: H^1(T) \times L^2(F) \to V_{T,F}^{H}$ which maps a given $v^H := (v , v|_{F}) \in H^1(T) \times L^2(F)$ to
\begin{align*}
\displaystyle
I_{T,F}^H v^H := (I_T^{CR} v , \Pi_F^0 v|_F).
\end{align*}
Furthermore, the global operator $I_h^H: X_0^H \to V_{dc,h,0}^H$ is defined as follows, for any $v^H := (v,v|_{\mathcal{F}_h}) \in X_0^H$,
\begin{align*}
\displaystyle
(I_h^H v^H)|_{T} = I_{T,F}^H (v^H|_T) \quad \forall T \in \mathbb{T}_h, \quad F \in \mathcal{F}_T.
\end{align*}

\begin{lem}[Best approximation] \label{best=approx}
We assume that $\Omega$ is convex. Let $u \in V_* := H_0^1(\Omega) \cap H^2(\Omega)$ be the solution of \eqref{poisson_weak} and $u^H := (u,u|_{\mathcal{F}_h}) \in X_0^H$. Then,
\begin{align}
\displaystyle
\inf_{v_h^H \in V_{dc,h,0}^{H}} |u^H - v_h^H|_{hwop(1)}
&\leq c \left(  \sum_{i=1}^d  \sum_{T \in \mathbb{T}_h} h_i^2 \left \| \frac{\partial}{\partial r_i} \nabla u \right \|_{L^2(T)^d}^2 \right)^{\frac{1}{2}}. \label{hwop=10}
\end{align}
\end{lem}

\begin{pf*}
Let $v \in H_0^1(\Omega)$ and $v^H := (v,v|_{\mathcal{F}_h}) \in X_0^H$. From the definitions of CR interpolation and $L^2$-projection, it holds that for any $T \in \mathbb{T}_h$ and $F \in \mathcal{F}_T$,
\begin{align}
\displaystyle
&\Pi_F^0 \left\{ (I_T^{CR} v - v) - (\Pi_F^0 v|_{\mathcal{F}_h} - v|_{\mathcal{F}_h}) \right\} \notag\\
&\quad = \frac{1}{|F|} \int_F \left\{ (I_T^{CR} v - v) - (\Pi_F^0 v|_{\mathcal{F}_h} - v|_{\mathcal{F}_h}) \right\} ds = 0. \label{hwop=11}
\end{align}
Therefore, using \eqref{CR4} and \eqref{hwop=11}, we obtain
\begin{align}
\displaystyle
\inf_{v_h^H \in V_{dc,h,0}^{H}} |u^H - v_h^H|_{hwop(1)}
&\leq |u^H - I_h^H u^H|_{hwop(1)} = |u - I_h^{CR} u|_{H^1(\mathbb{T}_h)} \notag \\
&\leq c \left(  \sum_{i=1}^d  \sum_{T \in \mathbb{T}_h} h_i^2 \left \| \frac{\partial}{\partial r_i} \nabla u \right \|_{L^2(T)^d}^2 \right)^{\frac{1}{2}}, \label{hwop=12new}
\end{align}
which is the target inequality \eqref{hwop=10}.
\qed
\end{pf*}

The essential part for error estimates is the consistency error term \eqref{hwop=9}.

\begin{lem}[Asymptotic consistency] \label{asy=con}
We assume that $\Omega$ is convex. Let $u \in V_*$ be the solution to \eqref{poisson_weak} and $u^H := (u,u|_{\mathcal{F}_h}) \in X_0^H$. Let $\{ \mathbb{T}_h\}$ be a family of conformal meshes with semi-regular properties (Assumption \ref{neogeo=assume}). Let $T \in \mathbb{T}_h$ be the element with Conditions \ref{cond1} or \ref{cond2} satisfying (Type \roman{sone}) in Section \ref{two=step} when $d=3$. Then,
\begin{align}
\displaystyle
E_h(u^H)
&\leq c \left\{ \left( \sum_{i=1}^d \sum_{T \in \mathbb{T}_h} h_i^2 \left \| \frac{\partial}{\partial r_i} \nabla u \right \|_{L^2(T)^d}^2 \right)^{\frac{1}{2}} + h \| \varDelta u \|_{L^2(\Omega)} \right\} \notag \\
&\quad + c \left ( h |u|_{H^1(\Omega)} + h^{\frac{3}{2}} |u|_{H^1(\Omega)}^{\frac{1}{2}} \| \varDelta u\|_{L^2(\Omega)}^{\frac{1}{2}}  \right). \label{hwop=12}
\end{align}
Furthermore, let $d=3$ and {$T \in \mathbb{T}_h$} be the element with Condition \ref{cond2} satisfying (Type \roman{stwo}) in Section \ref{two=step}. Then, it holds that
\begin{align}
\displaystyle
E_h(u^H)
&\leq c h \| \varDelta u \|_{L^2(\Omega)} + c \left ( h |u|_{H^1(\Omega)} + h^{\frac{3}{2}} |u|_{H^1(\Omega)}^{\frac{1}{2}} \| \varDelta u\|_{L^2(\Omega)}^{\frac{1}{2}}  \right). \label{hwop=13}
\end{align}
\end{lem}

\begin{pf*}
Let $u \in V_* := H_0^1(\Omega) \cap H^2(\Omega)$ and $u^H := (u,u|_{\mathcal{F}_h}) \in X_0^H$.  First, we have
\begin{align}
\displaystyle
 \div \mathcal{I}_{h}^{RT} \nabla u &= \Pi_h^{0} \div \nabla u =  \Pi_h^{0} \varDelta u, \label{rel=lem2}
\end{align}
{see \cite[Lemma 16.2]{ErnGue21a}.} Let $w_h^H = (w_h,\mu_h) \in V_{dc,h,0}^{H}$. Setting $w := \nabla u$ in \eqref{hwop=1} yields
\begin{align*}
\displaystyle
&a_h^{hwop}(u^H,w_h^H) - \ell_h(w_h^H) \\
&\quad =  \int_{\Omega} ( \nabla u- \mathcal{I}_{h}^{RT} \nabla u) \cdot \nabla_h w_{h} dx + \int_{\Omega} \left( \varDelta u -  \Pi_h^0 \varDelta u \right) w_{h}  dx  \\
&\quad \quad + \sum_{T \in \mathbb{T}_h} \sum_{F \in \mathcal{F}_T} n_{T,F} \int_{F} ( \nabla u \cdot n_F) \Pi_F^0 (w_{h} - \mu_h ) ds \\
&\quad =: I_1 + I_2 + I_3.
\end{align*}

Let $T \in \mathbb{T}_h$ be the element with Conditions \ref{cond1} or \ref{cond2} satisfying (Type \roman{sone}) in Section \ref{two=step} when $d=3$. Using the H\"older inequality, the Cauchy--Schwarz inequality, and the RT interpolation error \eqref{RT5}, the term $I_1$ is estimated as
\begin{align*}
\displaystyle
|I_1|
&\leq c  \sum_{T \in \mathbb{T}_h} \| \nabla u - \mathcal{I}_{h}^{RT} \nabla u \|_{L^2(T)} | w_{h} |_{H^1(T)} \\
&\leq c \sum_{T \in \mathbb{T}_h}\left(  \sum_{i=1}^d h_i \left \| \frac{\partial}{\partial r_i} \nabla u \right \|_{L^2(T)^d} + h_T \| \varDelta u \|_{L^2(T)} \right) | w_{h} |_{H^1(T)} \\
&\leq c \left\{ \left( \sum_{i=1}^d \sum_{T \in \mathbb{T}_h} h_i^2 \left \| \frac{\partial}{\partial r_i} \nabla u \right \|_{L^2(T)^d}^2 \right)^{\frac{1}{2}} + h \| \varDelta u \|_{L^2(\Omega)} \right\}  | w_{h}^H |_{hwop(1)}.
\end{align*}
Using the H\"older inequality, the Cauchy--Schwarz inequality, the stability of $\Pi_h^0$, and the estimate \eqref{L2ortho}, the term $I_2$ is estimated as
\begin{align*}
\displaystyle
|I_2|
&= \left|  \int_{\Omega} \left( \varDelta u -  \Pi_h^0 \varDelta u \right) \left(  w_{h} - \Pi_h^0 w_{h}  \right) dx \right| \\
&\leq \sum_{T \in \mathbb{T}_h} \| \varDelta u -  \Pi_h^0 \varDelta u \|_{L^2(T)}  \| w_{h} - \Pi_h^0 w_{h} \|_{L^2(T)} \\
&\leq c h \| \varDelta u \|_{L^2(\Omega)}  | w_{h}^H |_{hwop(1)}.
\end{align*}
Using inequality \eqref{hwop=2} with $\beta = 1$, the term $I_3$ is estimated as
\begin{align*}
\displaystyle
|I_3|
&\leq c \left(  h | u |_{H^1(\Omega)} + h^{ \frac{3}{2}} | u |_{H^1(\Omega)}^{\frac{1}{2}} \| \varDelta u \|_{L^2(\Omega)}^{\frac{1}{2}} \right) \notag \\
&\quad \quad \times  \left(  \sum_{T \in \mathbb{T}_h} \sum_{F \in \mathcal{F}_T} h^{- 2} \ell_{T,F}^{- 1} \left \|  \Pi_F^0 (w_{h} - \mu_h ) \right\|_{L^2(F)}^2 \right)^{\frac{1}{2}} \\
&\leq c \left(  h | u |_{H^1(\Omega)} + h^{ \frac{3}{2}} | u |_{H^1(\Omega)}^{\frac{1}{2}} \| \varDelta u \|_{L^2(\Omega)}^{\frac{1}{2}} \right)   | w_{h}^H |_{hwop(1)}.
\end{align*}

Furthermore, let $d=3$ and {$T \in \mathbb{T}_h$} be the element with Condition \ref{cond2} satisfying (Type \roman{stwo}) in Section \ref{two=step}. Then, it holds that
Using the H\"older inequality, the Cauchy--Schwarz inequality, and the RT interpolation error \eqref{RT6}, the term $I_1$ is estimated as
\begin{align*}
\displaystyle
|I_1|
&\leq c h  |  u |_{H^2(\Omega)}  | w_{h}^H |_{hwop(1)}, 
\end{align*}
which leads to the target inequality \eqref{hwop=13} with \eqref{poisson=reg}.
\qed
\end{pf*}

The following energy norm error estimates are the main results of this study.

\begin{thr} \label{energy=thr}
We assume that $\Omega$ is convex. Let $u \in V_*$ be the solution to \eqref{poisson_weak} and $u^H := (u,u|_{\mathcal{F}_h}) \in X_0^H$. Let $\{ \mathbb{T}_h\}$ be a family of conformal meshes with semi-regular properties (Assumption \ref{neogeo=assume}). Let $u_h^H := (u_h^{hwop}, \lambda_h^{hwop}) \ \in V_{dc,h,0}^{H}$ be the solution of \eqref{hwop=4}. Let $T \in \mathbb{T}_h$ be the element with Conditions \ref{cond1} or \ref{cond2} satisfying (Type \roman{sone}) in Section \ref{two=step} when $d=3$. Then,
\begin{align}
\displaystyle
|u^H - u_h^{H}|_{hwop(1)}
&\leq c \left\{ \left( \sum_{i=1}^d \sum_{T \in \mathbb{T}_h} h_i^2 \left \| \frac{\partial}{\partial r_i} \nabla u \right \|_{L^2(T)^d}^2 \right)^{\frac{1}{2}} + h \| \varDelta u \|_{L^2(\Omega)} \right\} \notag \\
&\quad + c \left ( h |u|_{H^1(\Omega)} + h^{\frac{3}{2}} |u|_{H^1(\Omega)}^{\frac{1}{2}} \| \varDelta u\|_{L^2(\Omega)}^{\frac{1}{2}}  \right). \label{hwop=14}
\end{align}
Furthermore, let $d=3$ and {$T \in \mathbb{T}_h$} be the element with Condition \ref{cond2} satisfying (Type \roman{stwo}) in Section \ref{two=step}. Then, it holds that
\begin{align}
\displaystyle
|u^H - u_h^{H}|_{hwop(1)}
&\leq c h \| \varDelta u \|_{L^2(\Omega)} + c \left ( h |u|_{H^1(\Omega)} + h^{\frac{3}{2}} |u|_{H^1(\Omega)}^{\frac{1}{2}} \| \varDelta u\|_{L^2(\Omega)}^{\frac{1}{2}}  \right). \label{hwop=15}
\end{align}
\end{thr}

\begin{pf*}
From Lemmata \ref{second=starang}, \ref{best=approx} and \ref{asy=con}, the intended estimates are obtained.
\qed
\end{pf*}

\subsection{$L^2$ norm error estimate}
This section presents the $L^2$ error estimate for the HWOPSIP method.

\begin{thr} \label{L2=thr}
We assume that $\Omega$ is convex. Let $u \in V_*$ be the solution to \eqref{poisson_weak} and $u^H := (u,u|_{\mathcal{F}_h}) \in X_0^H$. Let $\{ \mathbb{T}_h\}$ be a family of conformal meshes with semi-regular properties (Assumption \ref{neogeo=assume}). Let $u_h^H := (u_h^{hwop}, \lambda_h^{hwop}) \ \in V_{dc,h,0}^{H}$ be the solution of \eqref{hwop=4}. Let $T \in \mathbb{T}_h$ be the element with Conditions \ref{cond1} or \ref{cond2} satisfying (Type \roman{sone}) in Section \ref{two=step} when $d=3$. Then,
\begin{align}
\displaystyle
\|u - u_h^{hwop} \|_{L^2(\Omega)}
&\leq c h \left\{ \left( \sum_{i=1}^d \sum_{T \in \mathbb{T}_h} h_i^2 \left \| \frac{\partial}{\partial r_i} \nabla u \right \|_{L^2(T)^d}^2 \right)^{\frac{1}{2}} + h \| \varDelta u \|_{L^2(\Omega)} \right\} \notag \\
&\quad + c h  \left ( h |u|_{H^1(\Omega)} + h^{\frac{3}{2}} |u|_{H^1(\Omega)}^{\frac{1}{2}} \| \varDelta u\|_{L^2(\Omega)}^{\frac{1}{2}}  \right).  \label{hwop=17}
\end{align}
Furthermore, let $d=3$ and {$T \in \mathbb{T}_h$} be the element with Condition \ref{cond2} satisfying (Type \roman{stwo}) in Section \ref{two=step}. Then, it holds that
\begin{align}
\displaystyle
\|u - u_h^{hwop} \|_{L^2(\Omega)}
&\leq c h^2  \| \varDelta u \|_{L^2(\Omega)} \notag \\
&\quad + c h  \left ( h |u|_{H^1(\Omega)} + h^{\frac{3}{2}} |u|_{H^1(\Omega)}^{\frac{1}{2}} \| \varDelta u\|_{L^2(\Omega)}^{\frac{1}{2}}  \right).  \label{hwop=18}
\end{align}
\end{thr}
 
\begin{pf*}
We set $e:=  u - u_h^{hwop}$. Let $z \in V_*$ satisfy
\begin{align}
\displaystyle
- \varDelta z = e \quad \text{in $\Omega$}, \quad z = 0 \quad \text{on $\partial \Omega$}. \label{hwop=19}
\end{align}
We set $z^H := (z,z|_{\mathcal{F}_h})$. Note that $|z|_{H^1(\Omega)} \leq c \| e \|_{L^2(\Omega)}$ and $|z|_{H^2(\Omega)} \leq \| e \|_{L^2(\Omega)}$. Let $z_h^{H} := (z_h,  \eta_h) \in V_{dc,h,0}^{H}$ satisfy
\begin{align}
\displaystyle
a_h^{hwop}(\varphi_h^H , z_h^{H}) = \int_{\Omega} \varphi_h  e dx \quad \forall \varphi_h^H := (\varphi_h , \mu_h) \in V_{dc,h,0}^{H}. \label{hwop=20}
\end{align}
Then,
\begin{align*}
\displaystyle
\| e \|^2_{L^2(\Omega)}
&= \int_{\Omega} (u - u_h^{hwop}) e dx = a_h^{hwop}(u^H,z^H) - a_h^{hwop}(u_h^H , z_h^{H}) \\
&\hspace{-1cm} = a_{h}^{hwop}(u^H -  u_h^{H} , z^H - z_h^{H} ) + a_{h}^{hwop}( u^H -  u_h^{H} ,  z_h^{H}) + a_{h}^{hwop}( u_h^{H} ,  z^H - z_h^{H} ) \notag \\
&\hspace{-1cm} = a_{h}^{hwop}(u^H -  u_h^{H} , z^H - z_h^{H} ) \notag \\
&\hspace{-1cm} \quad + a_{h}^{hwop}( u^H -  u_h^{H} ,  z_h^{H} - I_{h}^{H} z^H) +   a_{h}^{hwop}( u^H -  u_h^{H} ,   I_{h}^{H} z^H) \notag\\
&\hspace{-1cm} \quad + a_{h}^{hwop}( u_h^{H} - I_{h}^{H} u^H ,  z^H - z_h^{H} ) +a_{h}^{hwop}( I_{h}^{H} u^H ,  z^H - z_h^{H} ) \\
&\hspace{-1cm} =: J_1 + J_2 + J_3 + J_4 + J_5.
\end{align*}

Let $T \in \mathbb{T}_h$ be the element with Conditions \ref{cond1} or \ref{cond2} satisfying (Type \roman{sone}) in Section \ref{two=step} when $d=3$.

Using \eqref{hwop=14}, $J_1$ can be estimated as
\begin{align*}
\displaystyle
|J_1| 
&\leq c |u^H - u_h^{H}|_{hwop(1)} |z^H - z_h^{H}|_{hwop(1)} \\
&\leq c h \| e \|_{L^2(\Omega)}  \left\{ \left( \sum_{i=1}^d \sum_{T \in \mathbb{T}_h} h_i^2 \left \| \frac{\partial}{\partial r_i} \nabla u \right \|_{L^2(T)^d}^2 \right)^{\frac{1}{2}} + h \| \varDelta u \|_{L^2(\Omega)} \right\} \notag \\
&\quad + c h \| e \|_{L^2(\Omega)} \left ( h |u|_{H^1(\Omega)} + h^{\frac{3}{2}} |u|_{H^1(\Omega)}^{\frac{1}{2}} \| \varDelta u\|_{L^2(\Omega)}^{\frac{1}{2}}  \right).
\end{align*}
Using \eqref{hwop=12new} and \eqref{hwop=14}, $J_2$ can be estimated as
\begin{align*}
\displaystyle
|J_2|
&\leq c | u^H- u_h^{H} |_{hwop(1)} \left(  | z_h^{H} - z^H |_{hwop(1)} + | z^H - I_h^{H} z^H|_{hwop(1)} \right) \\
&\leq c h \| e \|_{L^2(\Omega)}  \left\{ \left( \sum_{i=1}^d \sum_{T \in \mathbb{T}_h} h_i^2 \left \| \frac{\partial}{\partial r_i} \nabla u \right \|_{L^2(T)^d}^2 \right)^{\frac{1}{2}} + h \| \varDelta u \|_{L^2(\Omega)} \right\} \notag \\
&\quad + c h \| e \|_{L^2(\Omega)} \left ( h |u|_{H^1(\Omega)} + h^{\frac{3}{2}} |u|_{H^1(\Omega)}^{\frac{1}{2}} \| \varDelta u\|_{L^2(\Omega)}^{\frac{1}{2}}  \right).
\end{align*}
Using an analogous argument:
\begin{align*}
\displaystyle
|J_4|
&\leq c  \left(  | u_h^{H} - u^H |_{hwop(1)} + | u^H - I_h^{H} u^H|_{hwop(1)} \right)  | z^H- z_h^{H} |_{hwop(1)} \\
&\leq c h \| e \|_{L^2(\Omega)}  \left\{ \left( \sum_{i=1}^d \sum_{T \in \mathbb{T}_h} h_i^2 \left \| \frac{\partial}{\partial r_i} \nabla u \right \|_{L^2(T)^d}^2 \right)^{\frac{1}{2}} + h \| \varDelta u \|_{L^2(\Omega)} \right\} \notag \\
&\quad + c h \| e \|_{L^2(\Omega)} \left ( h |u|_{H^1(\Omega)} + h^{\frac{3}{2}} |u|_{H^1(\Omega)}^{\frac{1}{2}} \| \varDelta u\|_{L^2(\Omega)}^{\frac{1}{2}}  \right).
\end{align*}
Using $\mathcal{I}_h^{RT} \nabla u \cdot n_F \in \mathbb{P}^0(F)$ for any $F \in \mathcal{F}_h$ and the definition of CR interpolation,
\begin{align}
\displaystyle
&\int_{\Omega} \left( \mathcal{I}_h^{RT} \nabla u \cdot \nabla_h (I_{h}^{CR} z - z ) + \div \mathcal{I}_h^{RT} (\nabla u)  (I_{h}^{CR} z - z ) \right) dx \notag\\
&\quad =  \sum_{T \in \mathbb{T}_h} \sum_{F \in \mathcal{F}_T} n_{T,F} \int_{F} ( \mathcal{I}_h^{RT} \nabla u \cdot n_F) (I_{h}^{CR} z - z ) ds = 0. \label{hwop=21}
\end{align}
Furthermore, for any $T \in \mathbb{T}_h$ and $F \in \mathcal{F}_T$, we have
\begin{align}
\displaystyle
\Pi_F^0 (I_T^{CR} z|_F - \Pi_F^0 z|_F) = 0. \label{hwop=22}
\end{align}
Thus, using \eqref{rel=lem2}, \eqref{hwop=21}, and \eqref{hwop=22}, we have
\begin{align*}
\displaystyle
J_3
&=  a_{h}^{hwop}( u^H ,   I_{h}^{H} z^H) -  a_{h}^{hwop}( u_h^{H} ,   I_{h}^{H} z^H)\\
&= \int_{\Omega} \nabla u \cdot \nabla_h ( I_{h}^{CR} z - z ) dx + \int_{\Omega} \varDelta u (I_{h}^{CR} z - z )dx \\
&=  \int_{\Omega} (\nabla u -  \mathcal{I}_h^{RT} \nabla u  ) \cdot \nabla_h ( I_{h}^{CR} z - z ) dx + \int_{\Omega} ( \varDelta u - \Pi_h^0 \varDelta u ) (I_{h}^{CR} z - z )dx.
\end{align*}
Using \eqref{CR4}, \eqref{L2=CR4}, \eqref{RT5}, and the stability of the $L^2$-projection, we have
\begin{align*}
\displaystyle
|J_3|
&\leq  c h  \| e \|_{L^2(\Omega)} \left\{ \left( \sum_{i=1}^d \sum_{T \in \mathbb{T}_h} h_i^2 \left \| \frac{\partial}{\partial r_i} \nabla u \right \|_{L^2(T)^d}^2 \right)^{\frac{1}{2}} + h \| \varDelta u \|_{L^2(\Omega)} \right\} \\
&\quad + c h^2  \| e \|_{L^2(\Omega)} \| \varDelta u \|_{L^2(\Omega)}.
\end{align*}
Using an analogous argument:
\begin{align*}
\displaystyle
|J_5| 
&= \left| \int_{\Omega} (\nabla z -  \mathcal{I}_h^{RT} \nabla z  ) \cdot \nabla_h ( I_{h}^{CR} u - u ) dx + \int_{\Omega} ( \varDelta z - \Pi_h^0 \varDelta z ) (I_{h}^{CR} u - u )dx \right| \\
&\leq  c h  \| e \|_{L^2(\Omega)}  \left(  \sum_{i=1}^d  \sum_{T \in \mathbb{T}_h} h_i^2 \left \| \frac{\partial}{\partial r_i} \nabla u \right \|_{L^2(T)^d}^2 \right)^{\frac{1}{2}} + c h^2  \| e \|_{L^2(\Omega)} \| \varDelta u \|_{L^2(\Omega)}.
\end{align*}
Cobining the inequalities above yields the target inequality \eqref{hwop=17}.

Let $d=3$ and {$T \in \mathbb{T}_h$} be the element with Condition \ref{cond2} satisfying (Type \roman{stwo}) in Section \ref{two=step}. We use \eqref{hwop=15} instead of \eqref{hwop=14}, and \eqref{RT6} instead of \eqref{RT5} for the estimates of $J_1 - J_5$:
\begin{align*}
\displaystyle
|J_1|
&\leq c h^2  \| e \|_{L^2(\Omega)} \| \varDelta u \|_{L^2(\Omega)} + c h  \| e \|_{L^2(\Omega)} \left ( h |u|_{H^1(\Omega)} + h^{\frac{3}{2}} |u|_{H^1(\Omega)}^{\frac{1}{2}} \| \varDelta u\|_{L^2(\Omega)}^{\frac{1}{2}}  \right), \\
|J_2|
&\leq c h^2  \| e \|_{L^2(\Omega)} \| \varDelta u \|_{L^2(\Omega)} + c h  \| e \|_{L^2(\Omega)} \left ( h |u|_{H^1(\Omega)} + h^{\frac{3}{2}} |u|_{H^1(\Omega)}^{\frac{1}{2}} \| \varDelta u\|_{L^2(\Omega)}^{\frac{1}{2}}  \right), \\
|J_3|
&\leq c h^2  \| e \|_{L^2(\Omega)} \| \varDelta u \|_{L^2(\Omega)},\\
|J_4|
&\leq c h^2  \| e \|_{L^2(\Omega)} \| \varDelta u \|_{L^2(\Omega)} + c h  \| e \|_{L^2(\Omega)} \left ( h |u|_{H^1(\Omega)} + h^{\frac{3}{2}} |u|_{H^1(\Omega)}^{\frac{1}{2}} \| \varDelta u\|_{L^2(\Omega)}^{\frac{1}{2}}  \right), \\
|J_5|
&\leq c h^2  \| e \|_{L^2(\Omega)} \| \varDelta u \|_{L^2(\Omega)}.
\end{align*}
Combining the inequalities above yields the target inequality \eqref{hwop=18}.
\qed
\end{pf*}

\section{Numerical experiments} \label{numerical=sec}
This section presents the results of numerical experiments. Let $d=2$ and $\Omega := (0,1)^2$. 

\subsection{Example 1}
The first example shows the results of the numerical tests on anisotropic meshes. Function $f$ of the Poisson equation \eqref{poisson_eq},
\begin{align*}
\displaystyle
-  \varDelta u = f \quad \text{in $\Omega$}, \quad u = 0 \quad \text{on $\partial \Omega$}, 
\end{align*}
is given such that the exact solution is
{
\begin{align*}
\displaystyle
u (x_1,x_2) := 64 x_1 (x_1 - 1) x_2(x_2 - 1).
\end{align*}
}

\begin{figure}[htbp]
  \begin{minipage}[b]{0.45\linewidth}
    \centering
    \includegraphics[keepaspectratio, scale=0.15]{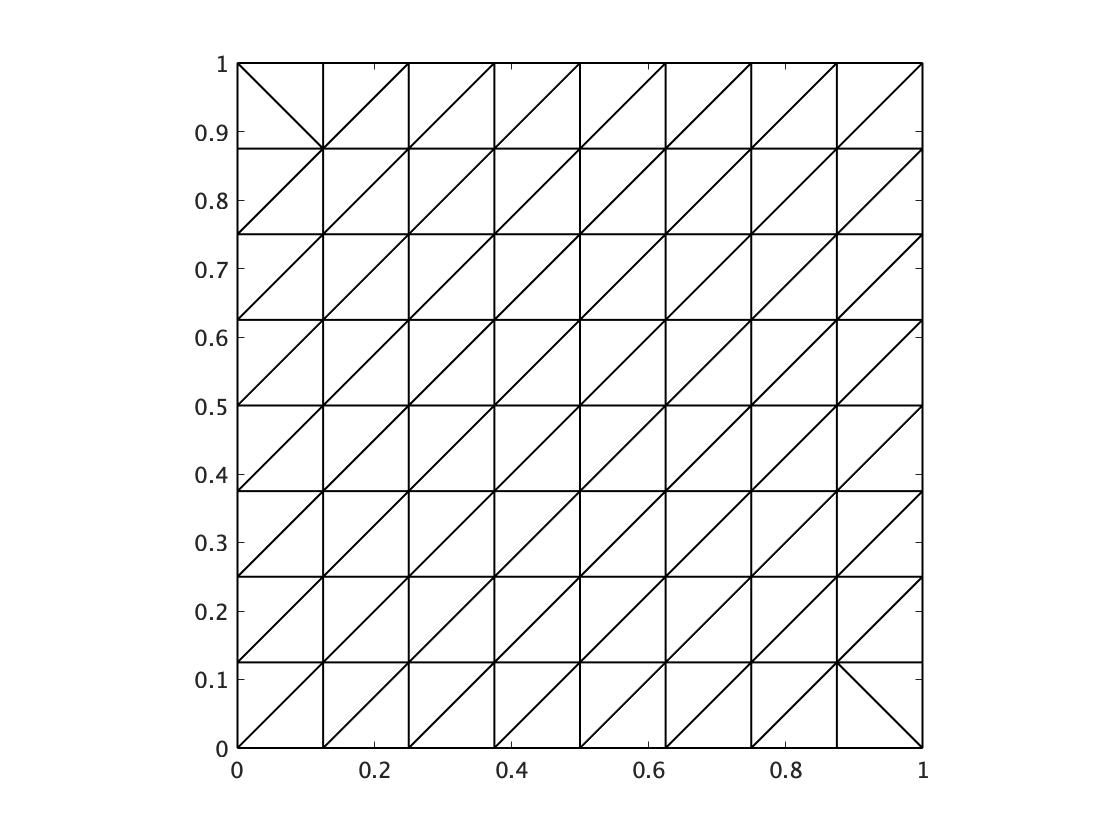}
    \caption{(\Roman{lone}) Standard mesh}
     \label{fig1}
  \end{minipage}
  \begin{minipage}[b]{0.45\linewidth}
    \centering
    \includegraphics[keepaspectratio, scale=0.15]{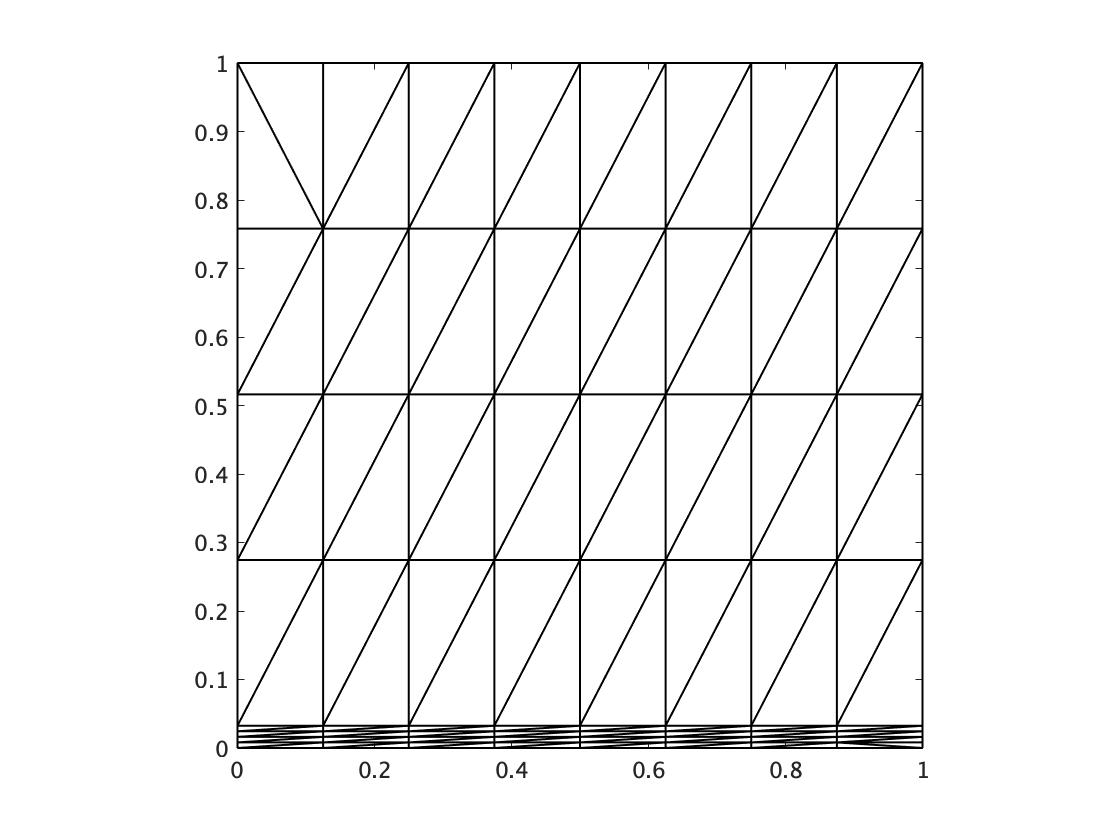}
    \caption{(\Roman{ltwo}) Shishkin mesh, $\delta = \frac{1}{128}$}
     \label{fig2}
  \end{minipage}
\end{figure}

\begin{figure}[htbp]
  \begin{minipage}[b]{0.45\linewidth}
    \centering
    \includegraphics[keepaspectratio, scale=0.15]{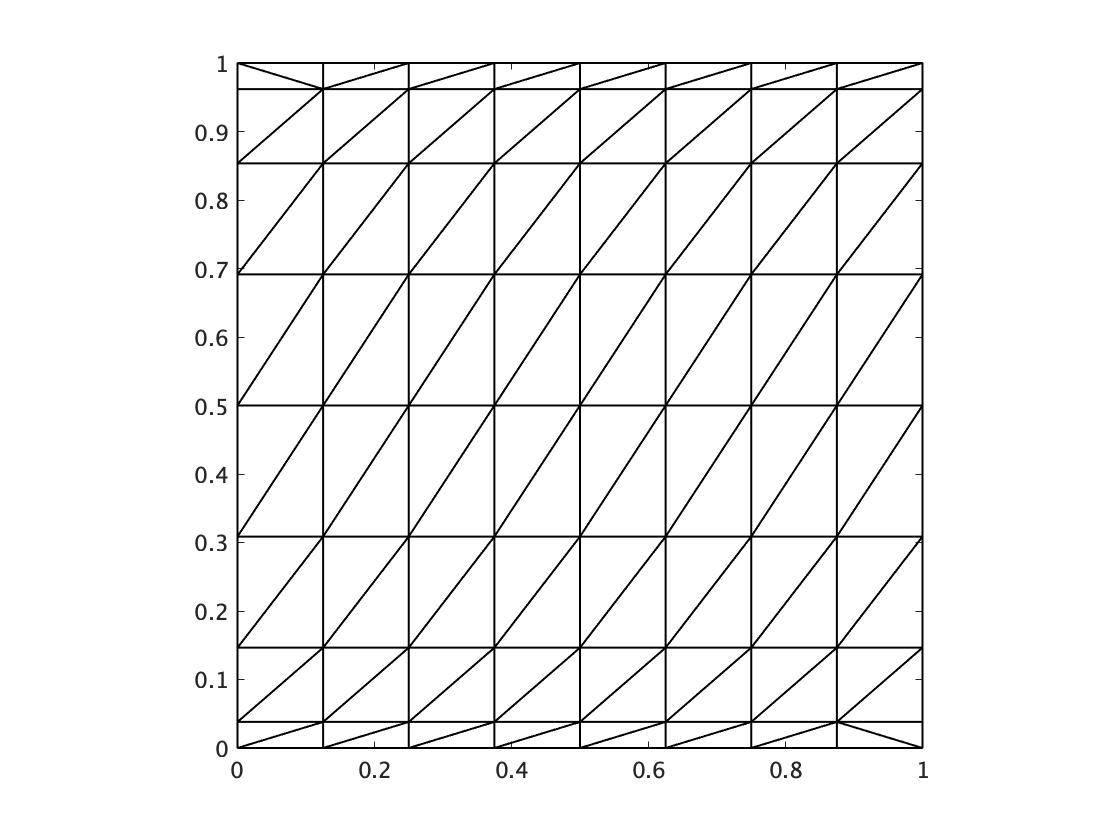}
    \caption{(\Roman{lthree}) Anisotropic mesh}
     \label{fig33}
  \end{minipage}
  \begin{minipage}[b]{0.45\linewidth}
    \centering
    \includegraphics[keepaspectratio, scale=0.15]{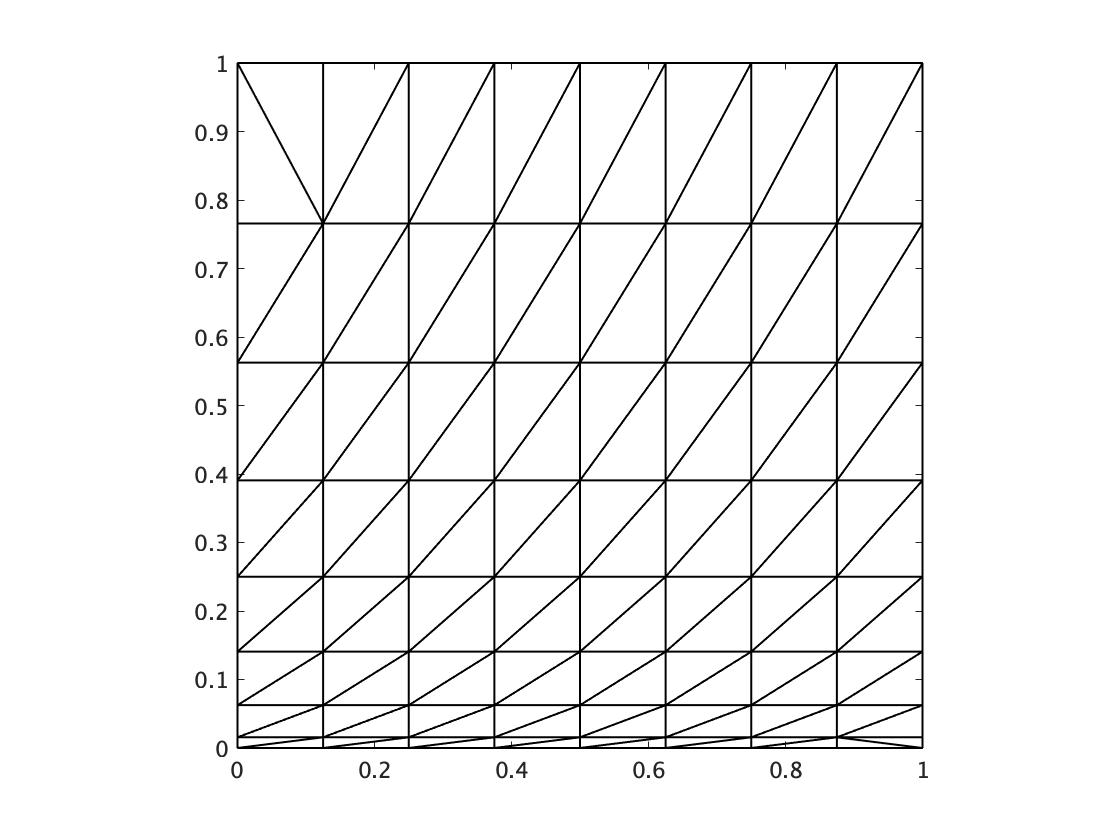}
    \caption{(\Roman{lfour}) {Graded mesh}}
     \label{fig44}
  \end{minipage}
\end{figure}

Let {$N \in \{ 32,64,128,256 \}$} be the division number of each side of the bottom and height edges of $\Omega$. Four types of mesh partitioning were considered. Let $(x_1^i, x_2^i)^T$ be the grip points of the triangulations $\mathbb{T}_h$ defined as follows: Let $i \in \mathbb{N}$.
\begin{description}
  \item[(\Roman{lone}) Standard mesh (Fig. \ref{fig1})] 
\begin{align*}
\displaystyle
x_1^i := \frac{i}{N}, \quad x_2^i := \frac{i}{N}, \quad  i \in \{0, \ldots , N \}.
\end{align*}
  \item[(\Roman{ltwo}) Shishkin mesh (Fig. \ref{fig2})]
\begin{align*}
\displaystyle
x_1^i &:= \frac{i}{N}, \quad  i \in \{0 , \ldots , N \}, \\
 x_2^i &:=
\begin{cases}
\tau \frac{2}{N} i, \quad  i \in \left\{0, \ldots , \frac{N}{2} \right \}, \\
\tau + (1 - \tau) \frac{2}{N} \left( i - \frac{N}{2} \right), \quad i \in \left\{ \frac{N}{2}+1 , \ldots , N \right\},
\end{cases}
\end{align*}
where $\tau := {2} \delta | \ln(N) |$ with $\delta = \frac{1}{128}$; see \cite[Section 2.1.2]{Lin10}.
\item[(\Roman{lthree}) Anisotropic mesh {from \cite{CheLiuQia10}} (Fig. \ref{fig33})]
\begin{align*}
\displaystyle
x_1^i := \frac{i}{N}, \quad x_2^i := \frac{1}{2}\left( 1 - \cos \left( \frac{i \pi}{N} \right) \right), \quad  i \in \{0, \ldots, N \}.
\end{align*}
 \item[(\Roman{lfour}) {Graded mesh} (Fig. \ref{fig44})]
 \begin{align*}
\displaystyle
x_1^i := \frac{i}{N}, \quad x_2^i := \left ( \frac{i}{N} \right)^{2}, \quad  i \in \{0, \ldots, N \}.
\end{align*}
\end{description}

Notably, a sequence with meshes (\Roman{lone}) or (\Roman{ltwo}) satisfies the shape-regularity condition, but a sequence with meshes (\Roman{lthree}) or (\Roman{lfour}) does not satisfy the shape-regularity condition; see \cite[Section 4.1]{Ish24b}.

If the exact solution $u$ is known, then the errors ${e_N := u - u_N}$ and ${e_{2N} := u - u_{2N}}$ are computed numerically for the {two division numbers $N$ and $2N$.} The convergence indicator $r$ is defined as
\begin{align*}
\displaystyle
r = \frac{1}{\log(2)} \log \left( \frac{\| {e_N} \|_X}{\| {e_{2N}} \|_X} \right).
\end{align*}
We compute the convergence order with respect to the norms defined by
\begin{align*}
\displaystyle
E_{H^1}^{(\text{Mesh No.})} &:= \frac{|u^H - {u_N^{H}} |_{hwop(1)} }{|u|_{H^1(\Omega)}}, \quad E_{L^2}^{(\text{Mesh No.})} :=  \frac{ \| u - {u_N} \|_{L^2(\Omega)}}{ \| u \|_{L^2(\Omega)}}.
\end{align*}
For the computation, we used the conjugate gradient method without preconditioners and the fifth-order quadrature for the computation of the right-hand side in \eqref{hwop=4}; see \cite[p. 85, Table 30.1]{ErnGue21b}. The notation $\# Np$ denotes the number of nodal points on $V_{dc,h,0}^{H}$, including the nodal points on the boundary. 

\begin{table}[htb]
\caption{Mesh (\Roman{lone})}
\centering
\begin{tabular}{l|l | l || l |l||l|l} \hline
$N$ & $\# Np$&  $h $ & $E_{H^1}^{(\Roman{lone})}$ & $r$ & $E_{L^2}^{(\Roman{lone})}$ & $r$  \\ \hline \hline
32& 9,280 & 4.42e-02  & 8.62073e-02   &   & 6.25004e-03    &   \\
64 &  36,992 &2.21e-02   &  4.31061e-02  & 1.00  &  1.56273e-03    & 2.00     \\
128 &  147,712 & 1.10e-02 & 2.15533e-02    &  1.00 &   3.90696e-04    &  2.00  \\
256 &590,336 &5.52e-03   & 1.07767e-02    &  1.00 &  9.76747e-05   &  2.00  \\
\hline
\end{tabular}
\label{mesh1}
\end{table}

\begin{table}[htb]
\caption{Mesh (\Roman{ltwo})}
\centering
\begin{tabular}{l|l | l || l |l||l|l} \hline
$N$ & $\# Np$&  $h $ & $E_{H^1}^{(\Roman{ltwo})}$ & $r$ & $E_{L^2}^{(\Roman{ltwo})}$ & $r$  \\ \hline \hline
32& 9,280 & 6.69e-02 & 1.31193e-01   &   & 1.45151e-02      &    \\
64 &  36,992 & 3.31e-02   &  6.50625e-02  & 1.01  &  3.57437e-03    & 2.02  \\
128 & 147,712 &1.64e-02  &   3.22586e-02 &  1.01 &     8.79693e-04  &  2.02  \\
256 & 590,336 & 8.13.e-03 & 1.59907e-02    & 1.01  & 2.16410e-04      &  2.02 \\
\hline
\end{tabular}
\label{mesh2}
\end{table}

\begin{table}[htb]
\caption{Mesh (\Roman{lthree})}
\centering
\begin{tabular}{l|l | l || l |l||l|l} \hline
$N$ & $\# Np$&  $h $ & $E_{H^1}^{(\Roman{lthree})}$ & $r$ & $E_{L^2}^{(\Roman{lthree})}$ & $r$  \\ \hline \hline
32&  9,280&  5.81e-02 & 1.12793e-01    &   & 1.09879e-02     &    \\
64 & 36,992 &  2.91e-02 &  5.64408e-02  &  1.00 &   2.75169e-03   & 2.00  \\
128 & 147,712 &1.45e-02   &  2.82259e-02  & 1.00  & 6.88219e-04     & 2.00 \\
256 & 590,336  & 7.27e-03  & 1.41137e-02  & 1.00  & 1.72073e-04     &  2.00  \\
\hline
\end{tabular}
\label{mesh3}
\end{table}

\begin{table}[htb]
\caption{Mesh (\Roman{lfour})}
\centering
\begin{tabular}{l|l | l || l |l||l|l} \hline
$N$ & $\# Np$&  $h $ & $E_{H^1}^{(\Roman{lfour})}$ & $r$ & $E_{L^2}^{(\Roman{lfour})}$ & $r$  \\ \hline \hline
32& 9,280 & 6.90e-02  & 1.30944e-01   &   &  1.53035e-02    &    \\
64 & 36,992 & 3.47e-02  &   6.58440e-02  &  0.99 &  3.87465e-03     &    1.98  \\
128 & 147,712 & 1.74e-02   & 3.30131e-02 & 1.00  &  9.74658e-04    &   1.99  \\
256 & 590,336  & 8.72e-03   & 1.65291e-02   & 1.00  &  2.44407e-04   &  2.00  \\
\hline
\end{tabular}
\label{mesh4}
\end{table}

Tables \ref{mesh1}--\ref{mesh4} present the numerical results of the HWOPSIP method for the standard and anisotropic meshes. We observed that all numerical experiments achieved the optimal convergence order. This demonstrates that the numerical solution converges to a sufficiently smooth exact solution even when the shape-regularity mesh condition is removed, but satisfies the semi-regular condition, which is also consistent with our theory.

\subsection{Example 2}
The second example demonstrates the numerical tests for a problem with boundary layers on meshes (\Roman{lone})–(\Roman{lfour}).  Function $f$ of the Poisson equation \eqref{poisson_eq} is given such that the exact solution is
{
\begin{align*}
\displaystyle
u (x_1,x_2) := 64 x_1 (x_1 - 1) x_2(x_2 - 1) \exp \left( - 128 x_2\right)
\end{align*}
}

\begin{table}[htb]
\caption{Mesh (\Roman{lone})}
\centering
\begin{tabular}{l|l | l || l |l||l|l} \hline
$N$ & $\# Np$&  $h $ & $E_{H^1}^{(\Roman{lone})}$ & $r$ & $E_{L^2}^{(\Roman{lone})}$ & $r$  \\ \hline \hline
32& 9,280 & 4.42e-02  &  1.09380  &   & 1.09891    &   \\
64 &  36,992 &2.21e-02   & 9.94176e-01   & 0.14  &  5.88132e-01   & 0.90     \\
128 &  147,712 & 1.10e-02 & 6.64606e-01  & 0.58   & 2.04815e-01   &  1.52   \\
256 &590,336 &5.52e-03   & 3.63751e-01  & 0.87  &  5.65848e-02    &  1.86   \\
\hline
\end{tabular}
\label{mesh1=2}
\end{table}

\begin{table}[htb]
\caption{Mesh (\Roman{ltwo})}
\centering
\begin{tabular}{l|l | l || l |l||l|l} \hline
$N$ & $\# Np$&  $h $ & $E_{H^1}^{(\Roman{ltwo})}$ & $r$ & $E_{L^2}^{(\Roman{ltwo})}$ & $r$  \\ \hline \hline
32& 9,280 & 6.39e-02 & 1.49440   &   &  1.72331   &   \\
64 &  36,992 & 3.14e-02   & 7.65292e-01  &  0.97 &  4.39920e-01   & 1.97  \\
128 & 147,712 &1.54e-02  &  3.88341e-01  & 0.98  &   1.11223e-01  &  1.98  \\
256 & 590,336 & 7.55e-03 &1.96481e-01    & 0.98  &    2.80265e-02  & 1.99  \\
\hline
\end{tabular}
\label{mesh2=2}
\end{table}

\begin{table}[htb]
\caption{Mesh (\Roman{lthree})}
\centering
\begin{tabular}{l|l | l || l |l||l|l} \hline
$N$ & $\# Np$&  $h $ & $E_{H^1}^{(\Roman{lthree})}$ & $r$ & $E_{L^2}^{(\Roman{lthree})}$ & $r$  \\ \hline \hline
32&  9,280&  5.81e-02 &  1.47146 &   &  1.68692    &    \\
64 & 36,992 &  2.91e-02 &  7.66982e-01   &  0.94 & 4.39689e-01      &  1.94\\
128 & 147,712 &1.45e-02   & 3.87770e-01   & 0.98  &  1.11162e-01    & 1.98 \\
256 & 590,336  & 7.27e-03  &1.94430e-01   & 1.00  &   2.78695e-02   & 2.00   \\
\hline
\end{tabular}
\label{mesh3=2}
\end{table}

\begin{table}[htb]
\caption{Mesh (\Roman{lfour})}
\centering
\begin{tabular}{l|l | l || l |l||l|l} \hline
$N$ & $\# Np$&  $h $ & $E_{H^1}^{(\Roman{lfour})}$ & $r$ & $E_{L^2}^{(\Roman{lfour})}$ & $r$  \\ \hline \hline
32& 9,280 & 6.90e-02  &  1.49372  &   &    1.72532   &    \\
64 & 36,992 & 3.47e-02  &  7.60015e-01  &  0.97 &    4.38991e-01 &   1.97   \\
128 & 147,712 & 1.74e-02   &  3.81716e-01 &  0.99 &     1.10245e-01&  1.99   \\
256 & 590,336  & 8.72e-03   & 1.91075e-01    & 1.00  &  2.75928e-02  &  2.00  \\
\hline
\end{tabular}
\label{mesh4=2}
\end{table}

Tables \ref{mesh1=2}--\ref{mesh4=2} present the numerical results of the HWOPSIP method for the boundary layer problem. We observe that the use of anisotropic meshes (\Roman{ltwo})–(\Roman{lfour}) near the bottom is likely to achieve the optimal convergence order, whereas Table \ref{mesh1=2} shows that the mesh needs to be sufficiently split to obtain the optimum convergence order on the standard mesh (\Roman{lone}). These results demonstrate the effectiveness of anisotropic meshes.

\subsection{Example 3}
The third example demonstrates numerical tests on meshes (\Roman{lone}) and (\Roman{lfour}) when the penalty parameters are changed. Function $f$ of the Poisson equation \eqref{poisson_eq} is given such that the exact solution is
\begin{align*}
\displaystyle
u (x_1,x_2) := 64 x_1 (x_1 - 1) x_2(x_2 - 1).
\end{align*}
We define the penalty parameters as
\begin{align*}
\displaystyle
\kappa_F^{(1)} := 10^{-2} \kappa_{F(1)}, \quad \kappa_F^{(2)} := 10^{2} \kappa_{F(1)}.
\end{align*}
We use the above parameters instead of $\kappa_{F(1)}$ in Scheme \eqref{hwop=4}.

\begin{table}[htb]
\caption{Mesh (\Roman{lone}), Case $\kappa_F^{(1)}$}
\centering
\begin{tabular}{l|l | l || l |l||l|l} \hline
$N$ & $\# Np$&  $h $ & $E_{H^1}^{(\Roman{lone})}$ & $r$ & $E_{L^2}^{(\Roman{lone})}$ & $r$  \\ \hline \hline
32& 9,280 & 4.42e-02  & 7.42966   &   & 5.66507e-01    &   \\
64 &  36,992 &2.21e-02   & 3.78302  & 0.97  & 1.44986e-01     &   1.97   \\
128 &  147,712 & 1.10e-02 & 1.90728  & 0.99   &  3.65684e-02   & 1.99   \\
256 &590,336 &5.52e-03   &9.56011e-01  &1.00   &   9.16498e-03  &   2.00 \\
\hline
\end{tabular}
\label{mesh1=3}
\end{table}

\begin{table}[htb]
\caption{Mesh (\Roman{lfour}), Case $\kappa_F^{(1)}$}
\centering
\begin{tabular}{l|l | l || l |l||l|l} \hline
$N$ & $\# Np$&  $h $ & $E_{H^1}^{(\Roman{lfour})}$ & $r$ & $E_{L^2}^{(\Roman{lfour})}$ & $r$  \\ \hline \hline
32& 9,280 & 6.90e-02  &  11.6940  &   & 1.40125     &    \\
64 & 36,992 & 3.47e-02  & 5.99162   & 0.96  &  3.65709e-01     & 1.94  \\
128 & 147,712 & 1.74e-02   & 3.04775   &  0.98 &   9.35448e-02   &  1.97 \\
256 & 590,336  & 8.72e-03   &   1.53430 & 0.99  &    2.35885e-02   & 1.99  \\
\hline
\end{tabular}
\label{mesh4=3}
\end{table}

\begin{table}[htb]
\caption{Mesh (\Roman{lone}), Case $\kappa_F^{(2)}$}
\centering
\begin{tabular}{l|l | l || l |l||l|l} \hline
$N$ & $\# Np$&  $h $ & $E_{H^1}^{(\Roman{lone})}$ & $r$ & $E_{L^2}^{(\Roman{lone})}$ & $r$  \\ \hline \hline
32& 9,280 & 4.42e-02  &3.96517e-02   &   &  1.18205e-03   &   \\
64 &  36,992 &2.21e-02   & 1.98342e-02  & 1.00  & 2.95815e-04   & 2.00     \\
128 &  147,712 & 1.10e-02 & 9.91797e-03 &  1.00 &  7.39710e-05  &  2.00  \\
256 &590,336 &5.52e-03   & 4.95908e-03  & 1.00  & 1.84938e-05   & 2.00   \\
\hline
\end{tabular}
\label{mesh1=4}
\end{table}

\begin{table}[htb]
\caption{Mesh (\Roman{lfour}), Case $\kappa_F^{(2)}$}
\centering
\begin{tabular}{l|l | l || l |l||l|l} \hline
$N$ & $\# Np$&  $h $ & $E_{H^1}^{(\Roman{lfour})}$ & $r$ & $E_{L^2}^{(\Roman{lfour})}$ & $r$  \\ \hline \hline
32& 9,280 & 6.90e-02  & 4.85273e-02   &   &  1.87536e-03    &    \\
64 & 36,992 & 3.47e-02  & 2.42901e-02  &  1.00 &   4.70667e-04    &  1.99 \\
128 & 147,712 & 1.74e-02   &  1.21480e-02  &  1.00  &  1.17789e-04    & 2.00  \\
256 & 590,336  & 8.72e-03   & 6.07435e-03  & 1.00  &    2.94683e-05   & 2.00  \\
\hline
\end{tabular}
\label{mesh4=4}
\end{table}

In our numerical tests (Tables \ref{mesh1=3}–\ref{mesh4=4}), an optimal convergence order was achieved without tuning any penalty parameters. However, changes in penalty parameters appear to affect the accuracy of the approximate solution. The use of parameter $\kappa_F^{(1)}$ increases the errors (Tables \ref{mesh1=3} and \ref{mesh4=3}), whereas the use of parameter $\kappa_F^{(2)}$ decreases the errors (Tables \ref{mesh1=4} and \ref{mesh4=4}). If the accuracy of the approximate solution for the original parameter $\kappa_{F(1)}$ is poor, changing the value of the parameter may reduce the errors.

\section{Concluding remarks} \label{sec=conc}
Here, we describe several topics for future research related to the results of this study.

\begin{enumerate}
 \item The HWOPSIP method proposed in this paper and the standard WOPSIP method may have similar advantages on anisotropic mesh partitions. A key advantage of the HWOPSIP method is that it does not require the tuning of any penalty parameter. Furthermore, the HWOPSIP method is expected to have the following good properties: (a) The appropriate Stokes element satisfies the inf-sup condition on anisotropic meshes; see \cite[Lemma 7]{Ish24}. (b) The error analysis of the method is studied on more general meshes than conformal meshes. This enables meshes with hanging nodes.
 \item A disadvantage of the HWOPSIP method is that it increases the number of conditions. The improvement of the ill-conditioning for the Poisson equation is stated in \cite{BreOweSun08} for the standard WOPSIP method. The improvement of the ill-conditioning owing to the over-penalty and the use of anisotropic meshes for the HWOPSIP method might be achievable as well.
 \item As described in Remark \ref{poin=rem1}, we impose that $\Omega$ is convex to obtain the stability of the discrete problem. The key idea is to deduce the discrete Poincar\'e and Sobolev inequalities on anisotropic meshes. However, to the best of our knowledge, the derivation of these inequalities remains an open question in non-conforming finite element methods on anisotropic meshes. 
 \item We used the discontinuous CR method to obtain error estimates. Extension to higher-order CR methods will be studied in the future.
\end{enumerate}

%
%


\begin{thebibliography}{}

\bibitem{AcoDur99}
Acosta, G., Dur\'an, R.G.: The maximum angle condition for mixed and nonconforming elements: Application to the Stokes equations. SIAM J. Numer. Anal. 37, 18--36 (1999)


\bibitem{BarBre14}
Barker, A. T., Brenner, S. C.: A mixed finite element method for the Stokes equations based on a weakly overpenalized symmetric interior penalty approach. J. Sci. Comput. 58, 290--307 (2014)

\bibitem{Bre03}
Brenner, S. C.: Poincar\'e--Friedrichs inequalities for piecewise $H^1$ functions. SIAM J. Numer. Anal. 41(1), 306--324 (2003)


\bibitem{Bre15}
Brenner, S. C.: Forty years of the Crouzeix--Raviart element. Numer. Methods Partial Differ. Equ. 31, 367--396 (2015) 

\bibitem{BreOweSun08}
Brenner, S. C., Owens, L., Sung, L. Y.: A weakly over-penalised symmetric interior penalty method. Electron. Trans. Numer. Anal. 30, 107--127 (2008)

\bibitem{BreOweSun12}
Brenner, S. C., Owens, L., Sung, L. Y.: Higher-order weakly overpenalized symmetric interior penalty methods. J. Comput. App. Math. 236, 2883--2894 (2012)

\bibitem{BreSco08}
Brenner, S. C., Scott, L. R.: The Mathematical Theory of Finite Element Methods. Springer-Verlag, New York (2008)

\bibitem{CheLiuQia10}
Chen, S., Liu, M., Qiao, Z.: Anisotropic nonconforming element for fourth-order elliptic singular perturbation problem. Int. J. Numer. Anal. Model. 7(4), 766--784 (2010)

\bibitem{Cocetal09}
Cockburn, B., Gopalakrichnan, J., Lazarov, R.: Unified hybridization of discontinuous Galerkin, mixed, and continuous Galerkin methods for second order elliptic problems. SIAM J. Numer. Anal. 47(2), 1319--1365 (2009)

\bibitem{CroRav73}
Crouzeix, M., Raviart, P.-A.: Conforming and nonconforming finite element methods for solving stationary Stokes equations, \Roman{lone}. RAIRO Anal. Num\'{e}r. 7, 33--76 (1973)

\bibitem{PieErn12}
Di Pietro, D. A., Ern, A.: Mathematical aspects of discontinuous Galerkin methods. Springer-Verlag (2012)

\bibitem{ErnGue04}
Ern, A., Guermond, J. L.: Theory and Practice of Finite Elements. Springer-Verlag, New York (2004)

\bibitem{ErnGue21a}
Ern, A., Guermond, J. L.: Finite Elements \Roman{lone}: Approximation and Interpolation. Springer-Verlag, New York (2021)

\bibitem{ErnGue21b}
Ern, A., Guermond, J. L.: Finite elements \Roman{ltwo}: Galerkin Approximation, Elliptic and Mixed PDEs. Springer-Verlag, New York (2021)

\bibitem{Gri11}
Grisvard, P.: Elliptical problems in Nonsmooth domains. SIAM (2011)

\bibitem{Ish21}
Ishizaka, H.: Anisotropic Raviart--Thomas interpolation error estimates using a new geometric parameter. Calcolo 59(4), (2022)

\bibitem{Ish24}
Ishizaka, H.: Anisotropic weakly over-penalised symmetric interior penalty method for the Stokes equation. Journal of Scientific Computing 100, (2024)

\bibitem{Ish24b}
Ishizaka, H.: Morley finite element analysis for fourth-order elliptic equations under a semi-regular mesh condition. Applications of Mathematics {69} (6), 769--805 (2024)

\bibitem{IshKobTsu21a}
Ishizaka, H., Kobayashi, K., Tsuchiya, T.: General theory of interpolation error estimates on anisotropic meshes. Jpn. J. Ind. Appl. Math. 38(1), 163--191 (2021)

\bibitem{IshKobSuzTsu21d}
Ishizaka, H., Kobayashi, K., Suzuki, R., Tsuchiya, T.: A new geometric condition equivalent to the maximum angle condition for tetrahedrons. Comput. Math. Appl. 99, 323--328 (2021)

\bibitem{IshKobTsu21c}
Ishizaka, H., Kobayashi, K., Tsuchiya, T.: Anisotropic interpolation error estimates using a new geometric parameter. Jpn. J. Ind. Appl. Math. 40(1), 475--512 (2023)

\bibitem{Lin10}
Lin\ss, T.: Layer-adapted meshes for reaction vector diffusion problems. Springer (2010)


\bibitem{Riv08}
Rivi\'ere, B.: Discontinuous Galerkin Methods for Solving Elliptic and Parabolic Equations. SIAM (2008)

\end{thebibliography}


\appendix
\section{Proof of the estimate \eqref{L2=CR4}} 
In this appendix, we provide the proof of  estimate \eqref{L2=CR4}:
\begin{align*}
\displaystyle
\|I_{T_j}^{CR} \varphi - \varphi \|_{L^2(T_j)}
&\leq c  \sum_{|\varepsilon| = 2} h^{\varepsilon} \left\| \partial_{r}^{\varepsilon} \varphi  \right\|_{L^{2}(T)} \quad \forall {\varphi} \in H^{2}({T}_j).
\end{align*}

\begin{pf*}
Let ${\varphi} \in H^{2}({T}_j)$ for $j \in \{ 1, \ldots , Ne\}$.  Using the scaling argument, we obtain
\begin{align}
\displaystyle
\| I_{T_j}^{CR} \varphi - \varphi \|_{L^2(T)}
\leq c |\det({A})|^{\frac{1}{2}} \| I_{\widehat{T}}^{CR} \hat{\varphi} - \hat{\varphi} \|_{L^2(\widehat{T})}. \label{th2=1}
\end{align}
For any $\hat{\eta} \in \mathbb{P}^{1}$, we have
\begin{align}
\displaystyle
\| I_{\widehat{T}}^{CR} \hat{\varphi} - \hat{\varphi} \|_{L^2(\widehat{T})}
&\leq \| I_{\widehat{T}}^{CR}  (\hat{\varphi} - \hat{\eta}) \|_{L^2(\widehat{T})}  + \| \hat{\eta} - \hat{\varphi} \|_{L^2(\widehat{T})},  \label{th2=2}
\end{align}
because $I_{\widehat{T}}^{CR} \hat{\eta} = \hat{\eta}$. Using the trace inequality in the reference element, 
\begin{align}
\displaystyle
\| I_{\widehat{T}}^{CR}  (\hat{\varphi} - \hat{\eta}) \|_{L^2(\widehat{T})} 
&\leq c \| \hat{\varphi} - \hat{\eta} \|_{H^1(\widehat{T})}. \label{th2=3}
\end{align}
From \eqref{th2=1}, \eqref{th2=2}, and \eqref{th2=3}, we obtain
\begin{align}
\displaystyle
\| I_{T_j}^{CR} \varphi - \varphi \|_{L^2(T)}
\leq c |\det({A})|^{\frac{1}{2}} \inf_{\hat{\eta} \in \mathbb{P}^{1}} \| \hat{\varphi} - \hat{\eta} \|_{H^1(\widehat{T})}. \label{th2=4}
\end{align}
From the Bramble--Hilbert lemma (\cite[Lemma 4.3.8]{BreSco08}), $\hat{\eta}_{\beta} \in \mathbb{P}^{1}$ exists such that for any $\hat{\varphi} \in H^{2}(\widehat{T})$,
\begin{align}
\displaystyle
| \hat{\varphi} - \hat{\eta}_{\beta} |_{H^{s}(\widehat{T})} \leq C^{BH}(\widehat{T}) |\hat{\varphi}|_{H^{2}(\widehat{T})}, \quad s=0,1. \label{th2=5}
\end{align}
Using the inequality in \cite[Lemma 6]{IshKobTsu21c} with $m=0$, we can estimate inequality \eqref{th2=5} as
\begin{align}
\displaystyle
|\hat{\varphi}|_{H^{2}(\widehat{T})}
\leq c |\det({{A}})|^{-\frac{1}{2}} \sum_{|\varepsilon| = 2} h^{\varepsilon} \left\| \partial_{r}^{\varepsilon} \varphi  \right\|_{L^{2}(T)}.  \label{th2=6}
\end{align}
By using \eqref{th2=4}, \eqref{th2=5} and \eqref{th2=6}, we can deduce the target inequality \eqref{L2=CR4}.
\qed
\end{pf*}

\end{document}